\documentclass[aos,MSNbibl,seceqn,dvips]{arximspdf}
\usepackage{subenv}
\usepackage{stfloats}
\usepackage{graphicx}

\doi{10.1214/12-AOS1018} %kopijuoti is PTS
\volume{40}
\issue{3}
\pubyear{2012}
\firstpage{1637}
\lastpage{1664}

\makeatletter

\fnbelowfloat

\newcommand{\Vvert}{|\!|\!|}

\newcommand{\order}{{\mathcal{O}}}
\newcommand{\Exs}{{\mathbb{E}}}

\newcommand{\defn}{:=}

\newcommand{\real}{{\mathbb{R}}}

\newcommand{\diag}{\operatorname{diag}}

\newcommand{\condind}{\perp\!\!\!\!\perp}
\newcommand{\betastar}{\beta^*}

\newcommand{\betahat}{{\widehat{\beta}}}
\newcommand{\nuhat}{{\widehat{\nu}}}
\newcommand{\ahat}{{\widehat{a}}}
\newcommand{\SigHat}{{\widehat{\Sigma}}}
\newcommand{\Thetahat}{{\widehat{\Theta}}}
\newcommand{\Gammahat}{{\widehat{\Gamma}}}
\newcommand{\Zbar}{{\overline{Z}}}

\newcommand{\betatil}{{\widetilde{\beta}}}
\newcommand{\Thetatil}{{\widetilde{\Theta}}}
\newcommand{\Ztil}{{\widetilde{Z}}}
\newcommand{\Mtil}{{\widetilde{M}}}

\newcommand{\rhobf}{{\bolds\rho}}
\newcommand{\rhohatbf}{{\widehat{\bolds\rho}}}

\newcommand{\numobs}{{n}}
\newcommand{\pdim}{{p}}

\newcommand{\E}{{\mathbb{E}}}

\newcommand{\kdim}{{k}}

\newcommand{\Qprob}{{\mathbb{Q}}}
\newcommand{\rhoprob}{{\rho}}
\newcommand{\CovW}{{\Sigma_w}}
\newcommand{\CovV}{{\Sigma_v}}
\newcommand{\CovX}{{\Sigma_x}}
\newcommand{\GamHat}{{\widehat{\Gamma}}}
\newcommand{\GamTil}{{\widetilde{\Gamma}}}

\newcommand{\GamLas}{{\GamHat_{\mathrm{Las}}}}

\newcommand{\GamMul}{{\GamHat_{\mathrm{mul}}}}

\newcommand{\GamMis}{{\GamHat_{\mathrm{mis}}}}

\newcommand{\GamAdd}{{\GamHat_{\mathrm{add}}}}

\newcommand{\Sset}{{S}}

\newcommand{\allow}{{\alpha_1}}
\newcommand{\alup}{{\alpha_2}}
\newcommand{\taulow}{{\tau}}
\newcommand{\tauup}{{\tau}}

\newcommand{\lammin}{{\lambda_{\mathrm{min}}}}
\newcommand{\lammax}{{\lambda_{\mathrm{max}}}}

\newcommand{\Sbar}{{\Sset^c}}

\newcommand{\sigeps}{{\sigma_\varepsilon}}

\newcommand{\spindex}{{k}}
\newcommand{\Loss}{{\mathcal{L}}}

\newcommand{\plaincon}{{c}}

\newcommand{\CONTRAC}{{\gamma}}

\newcommand{\rhomax}{{\rhoprob_{\mathrm{max}}}}
\newcommand{\rhoprobhat}{{\widehat{\rhoprob}}}

\newcommand{\regpar}{{\lambda_\numobs}}

\newcommand{\ANNOYCON}{{b_0}}

\newcommand{\hackdiv}{\mbox{\,${}:\hspace*{-8.2pt}\ominus{}$\,}}

\newproclaim{alg}{Algorithm}[section]
\newtheorem{theos}{Theorem}
\newtheorem{cors}{Corollary}

\newproclaim{exas}{Example}
\newproclaim{defns}{Definition}
\newproclaim{nota}{Notation}
\newproclaim{rem}{Remarks}

\makeatother

\begin{document}
\begin{frontmatter}

\title{High-dimensional regression with noisy and missing data:
Provable guarantees with nonconvexity}
\runtitle{High-dimensional noisy Lasso}

\begin{aug}
\author[A]{\fnms{Po-Ling} \snm{Loh}\corref{}\thanksref{t1,t2}\ead[label=e1]{ploh@berkeley.edu}}
\and
\author[A]{\fnms{Martin J.} \snm{Wainwright}\thanksref{t2}\ead[label=e2]{wainwrig@stat.berkeley.edu}}
\runauthor{P.-L. Loh and M. J. Wainwright}
\affiliation{University of California, Berkeley}
\address[A]{Department of Statistics\\
University of California, Berkeley \\
Berkeley, California 94720 \\
USA \\
\printead{e1}\\
\phantom{E-mail: }\printead*{e2}} %adresu isvedimo komanda gale!
\end{aug}

\thankstext{t1}{Supported in part by a Hertz Foundation Fellowship and
the Department of Defense (DoD) through an NDSEG Fellowship.}

\thankstext{t2}{Supported in part by NSF Grant DMS-09-07632 and Air
Force Office of Scientific Research Grant AFOSR-09NL184.}

% HISTORY:
\received{\smonth{9} \syear{2011}}
\revised{\smonth{5} \syear{2012}}

% ABSTRACT
%
\begin{abstract}
Although the standard formulations of prediction problems involve
fully-observed and noiseless data drawn in an i.i.d. manner, many
applications involve noisy and/or missing data, possibly involving
dependence, as well. We study these issues in the context of
high-dimensional sparse linear regression, and propose novel
estimators for the cases of noisy, missing and/or dependent data.
Many standard approaches to noisy or missing data, such as those using
the EM algorithm, lead to optimization problems that are inherently
nonconvex, and it is difficult to establish theoretical guarantees on
practical algorithms. While our approach also involves optimizing
nonconvex programs, we are able to both analyze the statistical error
associated with any global optimum, and more surprisingly, to prove
that a simple algorithm based on projected gradient descent will
converge in polynomial time to a small neighborhood of the set of all
global minimizers. On the statistical side, we provide nonasymptotic
bounds that hold with high probability for the cases of noisy,
missing and/or dependent data. On the computational side, we prove
that under the same types of conditions required for statistical
consistency, the projected gradient descent algorithm is guaranteed to
converge at a geometric rate to a near-global minimizer. We
illustrate these theoretical predictions with simulations, showing
close agreement with the predicted scalings.
\end{abstract}

% KEYWORDS
%
\begin{keyword}[class=AMS]
\kwd[Primary ]{62F12}
\kwd[; secondary ]{68W25}.
\end{keyword}
\begin{keyword}
\kwd{High-dimensional statistics}
\kwd{missing data}
\kwd{nonconvexity}
\kwd{regularization}
\kwd{sparse linear regression}
\kwd{$M$-estimation}.
\end{keyword}

\end{frontmatter}

%s1 #&#
\section{Introduction}

In standard formulations of prediction problems, it is assumed that
the covariates are fully-observed and sampled independently from some
underlying\vadjust{\goodbreak} distribution. However, these assumptions are not realistic
for many applications, in which covariates may be observed only
partially, observed subject to corruption or exhibit some type of
dependency. Consider the problem of modeling the voting behavior of
politicians: in this setting, votes may be missing due to abstentions,
and temporally dependent due to collusion or ``tit-for-tat'' behavior.
Similarly, surveys often suffer from the missing data problem, since
users fail to respond to all questions. Sensor network data also
tends to be both noisy due to measurement error, and partially missing
due to failures or drop-outs of sensors.

There are a variety of methods for dealing with noisy and/or missing
data, including various heuristic methods, as well as likelihood-based
methods involving the expectation--maximization (EM) algorithm (e.g.,
see the book~\cite{LitRub87} and references therein). A challenge in
this context is the possible nonconvexity of associated optimization
problems. For instance, in applications of EM, problems in which the
negative likelihood is a convex function often become nonconvex with
missing or noisy data. Consequently, although the EM algorithm will
converge to a local minimum, it is difficult to guarantee that the
local optimum is close to a global minimum.

In this paper, we study these issues in the context of
high-dimensional sparse linear regression---in particular, in the case
when the predictors or covariates are noisy, missing, and/or
dependent. Our main contribution is to develop and study simple
methods for handling these issues, and to prove theoretical results
about both the associated statistical error and the optimization
error. Like EM-based approaches, our estimators are based on solving
optimization problems that may be nonconvex; however, despite this
nonconvexity, we are still able to prove that a simple form of
projected gradient descent will produce an output that is
``sufficiently close''---as small as the statistical error---to any
global optimum. As a second result, we bound the statistical error,
showing that it has the same scaling as the minimax rates for the
classical cases of perfectly observed and independently sampled
covariates. In this way, we obtain estimators for noisy, missing,
and/or dependent data that have the same scaling behavior as the usual
fully-observed and independent case. The resulting estimators allow
us to solve the problem of high-dimensional Gaussian graphical model
selection with missing data.

There is a large body of work on the problem of corrupted covariates
or error-in-variables for regression problems (e.g., see the papers
and books~\cite{Hwa86,CarEtal95,ItuEtal99,XuYou07}, as well as
references therein). Much of the earlier theoretical work is
classical in nature, meaning that it requires that the sample size
$\numobs$ diverges with the dimension $\pdim$ fixed. Most relevant to
this paper is more recent work that has examined issues of corrupted
and/or missing data in the context of high-dimensional sparse linear
models, allowing for $\numobs\ll\pdim$. St\"{a}dler and
B\"{u}hlmann~\cite{StaBuh10} developed an EM-based method for sparse
inverse covariance matrix estimation in the missing data regime, and
used this result to derive an algorithm for sparse linear regression
with missing data. As mentioned above, however, it is difficult to
guarantee that EM will converge to a point close to a global optimum
of the likelihood, in contrast to the methods studied here. Rosenbaum
and Tsybakov~\cite{RosTsy10} studied the sparse linear model when the
covariates are corrupted by noise, and proposed a modified form of the
Dantzig selector (see the discussion
following our main results for a~detailed comparison to this past
work, and also to concurrent work~\cite{RosTsy11} by the same
authors). For the particular case of multiplicative noise, the type
of estimator that we consider here has been studied in past
work~\cite{XuYou07}; however, this theoretical analysis is of the
classical type, holding only for $\numobs\gg\pdim$, in contrast to
the high-dimensional models that are of interest here.

The remainder of this paper is organized as follows. We begin in
Section~\ref{SecBackground} with background and a precise description
of the problem. We then introduce the class of estimators we will
consider and the form of the projected gradient descent algorithm.
Section~\ref{SecMain} is devoted to a description of our main results,
including a pair of general theorems on the statistical and
optimization error, and then a series of corollaries applying our
results to the cases of noisy, missing, and dependent data. In
Section~\ref{SecSims}, we demonstrate simulations to confirm that our
methods work in practice, and verify the theoretically-predicted
scaling laws. Section~\ref{SecProofs} contains proofs of some of the main
results, with the remaining proofs contained in the supplementary
Appendix~\cite{LohWai11}.

\begin{nota*}
For a matrix $M$, we write $\|M\|_{\max} \defn\max_{i,j} |m_{ij}|$
to be the elementwise $\ell_\infty$-norm of $M$. Furthermore,
$\Vvert M \Vvert_{{1}}$ denotes the induced $\ell_1$-operator norm
(maximum absolute column sum) of $M$, and $\Vvert M \Vvert_{\mathrm
{op}}$ is the
spectral norm of $M$. We write $\kappa(M) \defn
\frac{\lammax(M)}{\lammin(M)}$, the condition number of~$M$. For
matrices $M_1, M_2$, we write $M_1 \odot M_2$ to denote the
componentwise Hadamard product, and write $M_1 \hackdiv M_2$ to
denote componentwise division. For functions $f(n)$ and $g(n)$, we
write $f(n) \precsim g(n)$ to mean that $f(n) \le c g(n)$ for
a universal constant $c \in(0, \infty)$, and similarly, $f(n)
\succsim g(n)$ when $f(n) \ge c' g(n)$ for some universal
constant $c' \in(0, \infty)$. Finally, we write $f(n) \asymp g(n)$
when $f(n) \precsim g(n)$ and $f(n) \succsim g(n)$ hold
simultaneously.
\end{nota*}

%s2 #&#
\section{Background and problem setup}
\label{SecBackground}

In this section, we provide background and a precise description of
the problem, and then motivate the class of estimators analyzed in
this paper. We then discuss a simple class of projected gradient
descent algorithms that can be used to obtain an estimator.

%s2.1 #&#
\subsection{Observation model and high-dimensional framework}

Suppose we observe a response variable $y_i \in\real$ linked to a
covariate vector $x_i \in\real^\pdim$ via the linear model
%
%e2.1 #&#
\begin{equation}
\label{LinMod} y_i = \bigl\langle x_i, \betastar\bigr\rangle+
\varepsilon_i \qquad\mbox{for $i = 1, 2, \ldots, \numobs$.}
\end{equation}
Here, the regression vector $\betastar\in\real^{\pdim}$ is unknown,
and $\varepsilon_i \in\real$ is observation noise, independent of
$x_i$. Rather than directly observing each $x_i \in\real^\pdim$, we
observe a vector $z_i \in\real^\pdim$ linked to $x_i$ via some
conditional distribution, that~is,
%
%e2.2 #&#
\begin{equation}
\label{EqnQprob} z_i \sim\Qprob(\cdot\mid x_i)
\qquad\mbox{for $i = 1, 2, \ldots, \numobs$.}
\end{equation}
This setup applies to various disturbances to the covariates,
including:
\begin{longlist}[(a)]
\item[(a)] \textit{Covariates with additive noise}: We observe $z_i =
x_i + w_i$, where $w_i \in\real^{\pdim}$ is a random vector
independent of $x_i$, say zero-mean with known covariance matrix~$\CovW$.
\item[(b)] \textit{Missing data}: For some fraction $\rhoprob\in
[0,1)$, we observe a random vector $z_i \in\real^{\pdim}$ such that
for each component $j$, we independently observe $z_{ij} = x_{ij}$
with probability $1-\rhoprob$, and $z_{ij} = \ast$ with
probability $\rhoprob$. We can also consider the case when the
entries in the $j$th column have a different probability
$\rhoprob_j$ of being missing.
\item[(c)] \textit{Covariates with multiplicative noise}: Generalizing
the missing data problem, suppose we observe $z_i = x_i \odot
u_i$, where $u_i \in\real^p$ is again a random vector independent
of $x_i$, and $\odot$ is the Hadamard product. The problem of
missing data is a special case of multiplicative noise, where all
$u_{ij}$'s are independent and $u_{ij} \sim
\operatorname{Bernoulli}(1-\rhoprob_j)$.
\end{longlist}
Our first set of results is deterministic, depending on specific
instantiations of the observations $\{(y_i, z_i)\}_{i=1}^n$. However,
we are also interested in results that hold with high probability when
the $x_i$'s and $z_i$'s are drawn at random. We consider both the case
when the $x_i$'s are drawn i.i.d. from a fixed distribution; and the
case of dependent covariates, when the $x_i$'s are generated according
to a stationary vector autoregressive (VAR) process.

We work within a high-dimensional framework that allows the number of
predictors $\pdim$ to grow and possibly exceed the sample size
$\numobs$. Of course, consistent estimation when $\numobs\ll\pdim$
is impossible unless the model is endowed with additional
structure---for instance, sparsity in the parameter vector
$\betastar$. Consequently, we study the class of models where
$\betastar$ has at most $\kdim$ nonzero parameters, where $\kdim$ is
also allowed to increase to infinity with $p$ and $n$. %\pldel{We
%implicitly assume the scaling $\|\betastar\|_2 = \order(1)$.}

%s2.2 #&#
\subsection{$M$-estimators for noisy and missing covariates}

In order to motivate the class of estimators we will consider, let us
begin by examining a~simple deterministic problem. Let $\CovX\succ
0$ be the covariance matrix of the covariates, and consider the
$\ell_1$-constrained quadratic program
%
%e2.3 #&#
\begin{equation}
\label{ExactProg} \betahat\in\mathop{\arg\min}_{\|\beta\|_1 \le R} \biggl\{ \frac{1}{2}
\beta^T \CovX\beta- \bigl\langle\CovX\betastar, \beta\bigr\rangle \biggr\}.
\end{equation}
As long as the constraint radius $R$ is at least $\|\betastar\|_1$,
the unique solution to this convex program is $\betahat= \betastar$.
Of course, this program is an idealization, since in practice we may
not know the covariance matrix $\CovX$, and we certainly do not know
$\CovX\betastar$---after all, $\betastar$ is the quantity we are
trying to estimate!\vadjust{\goodbreak}

Nonetheless, this idealization still provides useful intuition, as it
suggests various estimators based on the plug-in principle. Given a
set of samples, it is natural to form estimates of the quantities
$\CovX$ and $\CovX\betastar$, which we denote by $\GamHat\in
\real^{\pdim\times\pdim}$ and $\widehat{\gamma}\in
\real^\pdim$,
respectively, and to consider the modified program
%
%e2.4 #&#
\begin{equation}
\label{EqnGeneral} \betahat\in\mathop{\arg\min}_{\|\beta\|_1 \leq R} \biggl\{
\frac{1}{2} \beta^T \GamHat\beta- \langle\widehat{\gamma}, \beta
\rangle \biggr\},
\end{equation}
or alternatively, the regularized version
%
%e2.5 #&#
\begin{equation}
\label{EqnGeneralRegularized} \betahat\in\mathop{\arg\min}_{\beta\in\real^\pdim} \biggl\{
\frac{1}{2} \beta^T \GamHat\beta- \langle\widehat {\gamma},
\beta\rangle+ \regpar\|\beta\|_1 \biggr\},
\end{equation}
where $\regpar> 0$ is a user-defined regularization parameter. Note
that the two problems are equivalent by Lagrangian duality when the
objectives are convex, but not in the case of a nonconvex
objective. The Lasso~\cite{Tib96,CheEtal98} is a~special case of these
programs, obtained by setting
%
%e2.6 #&#
\begin{equation}
\label{EqnLassoChoice}
\GamLas\defn\frac{1}{\numobs} X^T X
\quad\mbox{and}\quad \widehat{\gamma}_{\mathrm{Las}} \defn\frac{1}{\numobs} X^T
y,
\end{equation}
where we have introduced the shorthand $y = (y_1, \ldots, y_\numobs)^T
\in\real^\numobs$, and $X \in\real^{\numobs\times\pdim}$, with
$x_i^T$ as its $i$th row. A simple calculation shows that
$(\GamLas, \widehat{\gamma}_{\mathrm{Las}})$ are unbiased
estimators of the pair $(\CovX,
\CovX\betastar)$. This unbiasedness and additional concentration
inequalities (to be described in the sequel) underlie the well-known
analysis of the Lasso in the high-dimensional regime.

In this paper, we focus on more general instantiations of the
programs~(\ref{EqnGeneral}) and (\ref{EqnGeneralRegularized}),
involving\vspace*{1pt} different choices of the pair $(\GamHat,
\widehat{\gamma})$ that are
adapted to the cases of noisy and/or missing data. Note that the
matrix $\GamLas$ is positive semidefinite, so the Lasso program is
convex. In sharp contrast, for the case of noisy or missing data, the
most natural choice of the matrix~$\GamHat$ is \textit{not positive
semidefinite}, hence the quadratic losses appearing in the
problems (\ref{EqnGeneral}) and (\ref{EqnGeneralRegularized}) are
\textit{nonconvex}. Furthermore, when $\GamHat$ has negative
eigenvalues, the objective in equation (\ref{EqnGeneralRegularized})
is unbounded from below. Hence, we make use of the following regularized
estimator:
%
%e2.7 #&#
\begin{equation}
\label{EqnGeneralSlight} \betahat\in
\mathop{\arg\min}_{\|\beta\|_1 \leq\ANNOYCON\sqrt{\kdim}} \biggl\{
\frac{1}{2} \beta^T \GamHat\beta- \langle \widehat{\gamma},
\beta\rangle+ \regpar \|\beta\|_1 \biggr\}
\end{equation}
for a suitable constant $\ANNOYCON$.

In the presence of nonconvexity, it is generally impossible to
provide a~polynomial-time algorithm that converges to a (near) global
optimum, due to the presence of local minima. Remarkably, we are able
to prove that this issue is not significant in our setting, and a
simple projected gradient descent algorithm applied to the
programs (\ref{EqnGeneral}) or (\ref{EqnGeneralSlight}) converges
with high probability to a vector extremely
close to any global optimum.

Let us illustrate these ideas with some
examples. Recall that $(\GamHat,
\widehat{\gamma})$ serve as unbiased
estimators for $(\Sigma_x, \Sigma_x \betastar)$.

%ex1 #&#
\begin{exas}[(Additive noise)]
\label{ExaNoisy}
Suppose we observe $Z = X+W$, where~$W$ is a random matrix independent
of $X$, with rows $w_i$ drawn i.i.d. from a~zero-mean distribution
with known covariance $\CovW$. We consider the pair
%
%e2.8 #&#
\begin{equation}
\label{EqnGamAdd} \GamAdd\defn\frac{1}{n} Z^TZ - \CovW
\quad\mbox{and}\quad \widehat{\gamma}_{\mathrm{add}} \defn\frac{1}{n}
Z^Ty.
\end{equation}
Note that when $\Sigma_w = 0$
(corresponding to the noiseless case), the estimators reduce to the
standard Lasso. However, when $\Sigma_w \neq0$, the matrix $\GamAdd$
is \textit{not} positive semidefinite in the high-dimensional regime ($n
\ll p$). Indeed, since the matrix $\frac{1}{n} Z^TZ$ has rank at most
$n$, the subtracted matrix $\Sigma_w$ may cause $\GamAdd$ to have a
large number of negative eigenvalues. For instance, if $\CovW=
\sigma^2_w I$ for $\sigma^2_w > 0$, then $\GamAdd$ has $\pdim-
\numobs$ eigenvalues equal to $-\sigma^2_w$.
\end{exas}

%%%%%%%%%%%%%%%%%%%%%%%%%%%%%%%%%%%%%%%%%%%%%%%%%%%%%%%%%%%%%%%%%%%%%%%%%%%%%%%%

%ex2 #&#
\begin{exas}[(Missing data)]
\label{ExaMiss}
We now consider the case where the entries of $X$ are missing at
random. Let us first describe an estimator for the special case where
each entry is missing at random, independently with some constant
probability $\rho\in[0,1)$. (In Example~\ref{ExaMul} to follow, we
will describe the extension to general missing probabilities.)
Consequently, we observe the matrix $Z \in\real^{\numobs\times
\pdim}$ with entries
\[
Z_{ij} = %
\cases{ X_{ij}, &\quad with probability $1-
\rho$,
\cr
0, &\quad otherwise.} %
\]
Given the observed matrix $Z \in\real^{\numobs\times\pdim}$, we use
%
%e2.9 #&#
\begin{equation}
\label{EqnSameGamMis} \GamMis\defn\frac{\Ztil^T \Ztil}{\numobs} - \rhoprob\diag \biggl(
\frac{\Ztil^T \Ztil}{\numobs} \biggr)
\quad\mbox{and}\quad \widehat{\gamma}_{\mathrm{mis}}\defn
\frac{1}{\numobs} \Ztil^T y,
\end{equation}
where $\Ztil_{ij} = Z_{ij}/(1-\rhoprob)$. It is easy to see that the
pair $(\GamMis, \widehat{\gamma}_{\mathrm{mis}})$ reduces to the
pair $(\GamLas,
\widehat{\gamma}_{\mathrm{Las}})$ for
the standard Lasso when $\rho= 0$, corresponding to no missing data.
In the more interesting case when $\rho\in(0,1)$, the matrix
$\frac{\Ztil^T \Ztil}{n}$ in equation (\ref{EqnSameGamMis}) has rank
at most $n$, so the subtracted diagonal matrix may cause the matrix
$\GamMis$ to have a large number of negative eigenvalues when $n \ll
p$. As a consequence, the matrix $\GamMis$ is not (in general)
positive semidefinite, so the associated quadratic function is not
convex.
\end{exas}

%%%%%%%%%%%%%%%%%%%%%%%%%%%%%%%%%%%%%%%%%%%%%%%%%%%%%%%%%%%%%%%%%%%%%%%%%%%%%%%%
%
%ex3 #&#
\begin{exas}[(Multiplicative noise)]
\label{ExaMul}
As a generalization of the previous example, we now consider the case
of multiplicative noise. In particular, suppose we observe the
quantity $Z = X \odot U$, where $U$ is a matrix of nonnegative noise
variables. In many applications, it is natural to assume that the
rows~$u_i$ of $U$ are drawn in an i.i.d. manner, say from some
distribution
in which both the vector $\Exs[u_1]$ and the matrix $\Exs[u_1 u_1^T]$
have strictly positive entries. This general family of multiplicative
noise models arises in various applications; we refer the reader to
the papers~\cite{Hwa86,CarEtal95,ItuEtal99,XuYou07} for more
discussion and examples. A natural\vspace*{1pt} choice of the pair $(\GamHat,
\widehat{\gamma})$ is given by the quantities
%
%e2.10 #&#
\begin{equation}
\label{EqnGamMul} \GamMul\defn\frac{1}{n} Z^TZ \hackdiv\E
\bigl(u_1u_1^T\bigr) \quad\mbox{and}\quad
\GamHat_{\mathrm{mul}}\defn\frac{1}{n} Z^Ty \hackdiv
\E(u_1),
\end{equation}
where $\hackdiv$ denotes elementwise division. A small calculation
shows that these are unbiased estimators of $\Sigma_x$ and
$\Sigma_x \betastar$, respectively. The estimators (\ref{EqnGamMul})
have been studied in past work~\cite{XuYou07}, but only under
classical scaling ($\numobs\gg\pdim$).

As a special case of the estimators (\ref{EqnGamMul}), suppose the
entries $u_{ij}$ of $U$ are independent $\operatorname{Bernoulli}(1-\rho_j)$ random
variables. Then the observed matrix $Z = X \odot U$ corresponds to a
missing-data matrix, where each element of the $j$th column
has probability $\rho_j$ of being missing. In this case, the
estimators~(\ref{EqnGamMul}) become
%
%e2.11 #&#
\begin{equation}
\label{EqnDefnGamMis} \GamMis= \frac{Z^TZ}{n} \hackdiv M
\quad\mbox{and}\quad \widehat{
\gamma}_{\mathrm{mis}}= \frac{1}{\numobs} Z^T y \hackdiv(\mathbf{1} -
\rhobf),
\end{equation}
where $M \defn\E(u_1u_1^T)$ satisfies
\[
M_{ij} = %
\cases{ (1-\rho_i) (1-
\rho_j), &\quad if $i \neq j$,
\cr
1-\rho_i, &\quad if $i
= j$,} %
\]
$\rhobf$ is the parameter vector containing the $\rho_j$'s, and
$\mathbf{1}$ is the vector of all 1's. In this way, we obtain a
generalization of the estimator discussed in Example~\ref{ExaMiss}.
\end{exas}

%s2.3 #&#
\subsection{Restricted eigenvalue conditions}

Given an estimate $\betahat$, there are various ways to assess its
closeness to $\betastar$. In this paper,\vspace*{1pt} we focus on the
$\ell_2$-norm $\|\betahat- \betastar\|_2$, as well as the closely
related $\ell_1$-norm $\|\betahat- \betastar\|_1$. When the
covariate matrix $X$ is fully observed (so that the Lasso can be
applied), it is now well understood that a sufficient condition for
$\ell_2$-recovery is that the matrix $\GamLas= \frac{1}{\numobs} X^T
X$ satisfy a certain type of restricted eigenvalue (RE) condition
(e.g.,~\cite{BicEtal08,GeeBuh09}). In this paper, we make use of the
following condition.

%de1 #&#
\begin{defns}[(Lower-RE condition)]
The matrix $\GamHat$ satisfies a lower restricted eigenvalue condition
with curvature $\allow> 0$ and tolerance \mbox{$\taulow(\numobs,
\pdim) > 0$} if
%
%e2.12 #&#
\begin{equation}
\label{EqnRSC} \theta^T \GamHat\theta\geq\allow\|\theta
\|_2^2 - \taulow(\numobs, \pdim) \|\theta
\|_1^2 \qquad\mbox{for all $\theta \in\real^\pdim$.}
\end{equation}
\end{defns}

It can be shown that when the Lasso matrix $\GamLas=
\frac{1}{\numobs} X^T X$ satisfies this RE condition (\ref{EqnRSC}),
the Lasso estimate has low $\ell_2$-error for any vector $\betastar$
supported on any subset of size at most $\spindex\lesssim
\frac{1}{\taulow(\numobs, \pdim)}$. In particular,
bound~(\ref{EqnRSC}) implies a sparse RE condition for all $\spindex$
of this magnitude, and conversely, Lemma 11 in\vadjust{\goodbreak}
the Appendix of~\cite{LohWai11} shows that a sparse RE condition implies
bound (\ref{EqnRSC}). In this paper, we work with
condition (\ref{EqnRSC}), since it is especially convenient for
analyzing optimization algorithms.

In the standard setting (with uncorrupted and fully observed design
matrices), it is known that for many choices of the design matrix $X$
(with rows having covariance $\Sigma$), the Lasso matrix $\GamLas$
will satisfy such an RE condition with high probability
(e.g.,~\cite{RasEtal10,RudZho11}) with $\allow= \frac{1}{2}
\lammin(\Sigma)$ and $\taulow(\numobs, \pdim) \asymp\frac{\log
\pdim}{\numobs}$. A significant portion of the analysis in this
paper is devoted to proving that different choices of $\GamHat$, such
as the matrices $\GamAdd$ and $\GamMis$ defined earlier, also satisfy
condition (\ref{EqnRSC}) with high probability. This fact is by
no means obvious, since as previously discussed, the matrices
$\GamAdd$ and $\GamMis$ generally have large numbers of negative
eigenvalues.

Finally, although such upper bounds are not necessary for statistical
consistency, our algorithmic results make use of the analogous upper
restricted eigenvalue condition, formalized in the following:
%
%de2 #&#
\begin{defns}[(Upper-RE condition)]
The matrix $\GamHat$ satisfies an upper restricted eigenvalue
condition with smoothness $\alup> 0$ and tolerance $\tauup(\numobs,
\pdim) > 0$ if
%
%e2.13 #&#
\begin{equation}
\label{EqnRSM} \theta^T \GamHat\theta\leq\alup\|\theta
\|_2^2 + \tauup(\numobs, \pdim) \|\theta
\|_1^2 \qquad\mbox{for all $\theta\in\real^\pdim$.}
\end{equation}
\end{defns}

In recent work on high-dimensional projected gradient
descent, Agarwal et al.~\cite{AgaEtal11} make use of a more general
form of the lower and upper bounds~(\ref{EqnRSC}) and (\ref{EqnRSM}),
applicable to nonquadratic losses as well, which are referred to as
the restricted strong convexity (RSC) and restricted smoothness (RSM)
conditions, respectively. For various class of random design
matrices, it can be shown that the Lasso matrix $\GamLas$ satisfies
the upper bound~(\ref{EqnRSM}) with $\alup= 2 \lammax(\CovX)$ and
$\tauup(\numobs, \pdim) \asymp\frac{\log\pdim}{\numobs}$; see
Raskutti et al.~\cite{RasEtal10} for the Gaussian case and Rudelson
and Zhou~\cite{RudZho11} for the sub-Gaussian setting. We
will establish similar scaling for our choices of $\GamHat$.
%for $(n,p,k)$ to be such that $\kdim\taulow(\numobs, \pdim)$ is
%bounded above by a small constant.}

%s2.4 #&#
\subsection{Gradient descent algorithms}

In addition to proving results about the global minima of the
(possibly nonconvex) programs (\ref{EqnGeneral})
and (\ref{EqnGeneralRegularized}), we are also interested in
polynomial-time procedures for approximating such optima. In
this paper, we analyze some simple algorithms for solving either the
constrained program (\ref{EqnGeneral}) or the Lagrangian
version (\ref{EqnGeneralSlight}). Note that the gradient of the
quadratic loss function takes the form $\nabla\Loss(\beta) =
\GamHat\beta- \widehat{\gamma}$. In application to
the constrained
version, the method of projected gradient descent generates a
sequence of iterates $\{\beta^t, t = 0, 1, 2, \ldots\}$ by the
recursion
%
%e2.14 #&#
\begin{equation}
\label{EqnProjGrad} %\beta^{t+1} & = \Pi\big(\beta^t - \frac{1}{\eta} (\GamHat\beta^t -
\beta^{t+1} = \mathop{\arg
\min}_{\|\beta\|_1 \leq R} \biggl\{\Loss\bigl(\beta^t\bigr) + \bigl\langle
\nabla\Loss\bigl(\beta^t\bigr), \beta- \beta^t \bigr
\rangle+ \frac{\eta}{2} \bigl\|\beta- \beta^t\bigr\|_2^2
\biggr\},
\end{equation}
where $\eta> 0$ is a stepsize parameter. Equivalently, this update
can be written as $\beta^{t+1} = \Pi (\beta^t - \frac{1}{\eta}
\nabla\Loss(\beta^t)  )$, where $\Pi$ denotes the\vadjust{\goodbreak}
$\ell_2$-projection onto the $\ell_1$-ball of radius $R$. This
projection can be computed rapidly in $\order(\pdim)$ time using a
procedure due to Duchi et al.~\cite{DucEtal08}. For the Lagrangian
update, we use a~slight variant of the projected gradient
update (\ref{EqnProjGrad}), namely
%
%e2.15 #&#
\begin{equation}
\label{EqnMirror} %\beta^{t+1} & = \Pi\big(\beta^t - \frac{1}{\eta} (\GamHat\beta^t -
\beta^{t+1} = \mathop{\arg
\min}_{\|\beta\|_1 \leq R} \biggl\{\Loss\bigl(\beta^t\bigr) + \bigl\langle
\nabla\Loss\bigl(\beta^t\bigr), \beta- \beta^t \bigr
\rangle+ \frac{\eta}{2} \bigl\|\beta- \beta^t\bigr\|_2^2
+ \regpar\|\beta\|_1 \biggr\}\hspace*{-35pt}
\end{equation}
with the only difference being the inclusion of the regularization
term. This update can also performed efficiently by performing two
projections onto the $\ell_1$-ball; see the paper~\cite{AgaEtal11} for
details.

When the objective function is convex (equivalently, $\GamHat$ is
positive semidefinite), the iterates (\ref{EqnProjGrad})
or (\ref{EqnMirror}) are guaranteed to converge to a global minimum of
the objective functions (\ref{EqnGeneral})
and (\ref{EqnGeneralSlight}), respectively. In our setting, the
matrix $\GamHat$ need not be positive semidefinite, so the best
generic guarantee is that the iterates converge to a local optimum.
However, our analysis shows that for the family of
programs (\ref{EqnGeneral}) or (\ref{EqnGeneralSlight}), under a
reasonable set of conditions satisfied by various statistical models,
the iterates actually converge to a point extremely close to any
global optimum in both $\ell_1$-norm and $\ell_2$-norm; see
Theorem~\ref{ThmOpt} to follow for a more detailed statement.

%s3 #&#
\section{Main results and consequences}
\label{SecMain}
We now state our main results and discuss their consequences for
noisy, missing, and dependent data.

%s3.1 #&#
\subsection{General results}

%We provide theoretical guarantees for both the constrained
%estimator \eqref{EqnGeneral} and the regularized variant
%%
%{\gamma}}}{\beta} + \regpar
%%
%for a constant $\ANNOYCON\ge\|\betastar\|_2$, which is a hybrid
%between the constrained \eqref{EqnGeneral} and
%regularized \eqref{EqnGeneralRegularized} programs, and convenient for
%technical purposes.
%}

We provide theoretical guarantees for both the constrained
estimator (\ref{EqnGeneral}) and the Lagrangian
version (\ref{EqnGeneralSlight}). Note that we obtain different
optimization problems as we vary the choice of the pair $(\GamHat,
\widehat{\gamma}) \in\real^{\pdim\times\pdim}
\times\real^\pdim$. We begin
by stating a pair of general results, applicable to any pair that
satisfies certain conditions. Our first result
(Theorem~\ref{ThmStat}) provides bounds on the statistical error,
namely the quantity $\|\betahat- \betastar\|_2$, as well as the
corresponding $\ell_1$-error, where $\betahat$ is any global optimum
of the programs (\ref{EqnGeneral}) or (\ref{EqnGeneralSlight}). Since
the problem may be nonconvex in general, it is not immediately
obvious that one can obtain a \textit{provably good} approximation to
any global optimum without resorting to costly search methods. In
order to assuage this concern, our second result
(Theorem~\ref{ThmOpt}) provides rigorous bounds on the optimization
error, namely the differences $\|\beta^t - \betahat\|_2$ and
$\|\beta^t - \betahat\|_1$ incurred by the iterate $\beta^t$ after
running $t$ rounds of the projected gradient descent
updates (\ref{EqnProjGrad}) or (\ref{EqnMirror}).

%s3.1.1 #&#
\subsubsection{Statistical error}

In controlling the statistical error, we assume that the matrix
$\GamHat$ satisfies a lower-RE condition with curvature $\allow$ and
tolerance $\taulow(\numobs, \pdim)$, as previously
defined (\ref{EqnRSC}). Recall that $\GamHat$ and $
\widehat{\gamma}$ serve as
surrogates to the deterministic quantities $\CovX\in\real^{\pdim
\times\pdim}$ and $\CovX\betastar\in\real^\pdim$, respectively.
Our results also involve a measure of deviation in these
surrogates. In particular, we assume that there is some function
$\varphi(\Qprob, \sigeps)$, depending on the two sources of noise in
our problem: the standard deviation $\sigeps$ of the observation noise
vector $\varepsilon$ from equation (\ref{LinMod}), and the conditional
distribution $\Qprob$ from equation (\ref{EqnQprob}) that links the
covariates $x_i$ to the observed versions $z_i$. With this notation,
we consider the deviation condition
%
%e3.1 #&#
\begin{equation}
\label{EqnDevFull} \bigl\|\widehat{\gamma}- \GamHat\betastar\bigr\|_\infty\le
\varphi(\Qprob, \sigeps) \sqrt{\frac{\log p}{n}}.
\end{equation}
To aid intuition, note that inequality (\ref{EqnDevFull}) holds
whenever the following two deviation conditions are satisfied:
%
%e3.2 #&#
\begin{eqnarray}
\label{EqnDev} \bigl\|\widehat{\gamma}- \CovX\betastar\bigr\|_\infty&\leq& \varphi(
\Qprob, \sigeps) \sqrt{\frac{\log\pdim}{\numobs}}
\quad\mbox{and}\nonumber\\[-8pt]\\[-8pt]
\bigl\| (\GamHat- \CovX)
\betastar\bigr\|_\infty&\leq&\varphi(\Qprob, \sigeps)
\sqrt{\frac{\log\pdim}{\numobs}}.\nonumber
\end{eqnarray}
The pair of inequalities (\ref{EqnDev}) clearly measures the
deviation of the estimators $(\GamHat, \widehat{\gamma
})$ from their
population versions, and they are sometimes easier to verify
theoretically. However, inequality (\ref{EqnDevFull}) may be used
directly to derive tighter bounds (e.g., in the additive noise
case). Indeed, the bounds established via
inequalities (\ref{EqnDev}) is not sharp in the limit of low noise
on the covariates, due to the second inequality. In the proofs of
our corollaries to follow, we will verify the deviation conditions for
various forms of noisy, missing, and dependent data, with the quantity
$\varphi(\Qprob, \sigeps)$ changing depending on the model. We have
the following result, which applies to any global optimum $\betahat$
of the regularized version (\ref{EqnGeneralSlight}) with
$\lambda_\numobs\ge4 \varphi(\Qprob, \sigeps)
\sqrt{\frac{\log\pdim}{\numobs}}$:
%
%th1 #&#
\begin{theos}[(Statistical error)]
\label{ThmStat}
Suppose the surrogates $(\GamHat, \widehat{\gamma})$
satisfy the deviation
bound (\ref{EqnDevFull}), and the matrix $\GamHat$ satisfies
the lower-RE condition~(\ref{EqnRSC}) with parameters $(\allow,
\taulow)$ such that
%
%e3.3 #&#
\begin{equation}
\label{EqnRSCConstants}
\sqrt{\kdim} \taulow(\numobs, \pdim) \leq\min \Biggl\{
\frac{\allow}{128 \sqrt{\kdim}}, \frac{\varphi(\Qprob,
\sigeps)}{\ANNOYCON} \sqrt{\frac{\log\pdim}{\numobs}} \Biggr\}.
\end{equation}
Then for any vector $\betastar$ with sparsity at most $\kdim$, there
is a universal positive constant $\plaincon_0$ such that any global
optimum $\betahat$ of the Lagrangian
program (\ref{EqnGeneralSlight}) with any $b_0 \ge
\|\betastar\|_2$ satisfies the bounds
%
%e3.4 #&#
\begin{subequation}
\label{EqnGenBound}
%
%e3.4a #&#
%e3.4b #&#
\begin{eqnarray}
\label{EqnEllTwoBound} \bigl\|\betahat- \betastar\bigr\|_2 & \leq &
\frac{\plaincon_0
\sqrt{\kdim}}{\allow} \max \Biggl\{ \varphi(\Qprob, \sigeps) \sqrt{
\frac{\log\pdim}{\numobs}}, \lambda_n \Biggr\}
\quad\mbox{and}\quad
\\
\label{EqnEllOneBound} \bigl\|\betahat- \betastar\bigr\|_1 & \leq &
\frac{8 \plaincon_0
\kdim}{\allow} \max \Biggl\{ \varphi(\Qprob, \sigeps) \sqrt{
\frac{\log\pdim}{\numobs}}, \lambda_n \Biggr\}.
\end{eqnarray}
\end{subequation}
\end{theos}

The same bounds (without $\regpar$) also apply to the constrained
program~(\ref{EqnGeneral}) with radius choice $R =
\|\betastar\|_1$.\vspace*{8pt}

\textit{Remarks.\quad}
To be clear, all the claims of Theorem~\ref{ThmStat} are
deterministic. Probabilistic conditions will enter when we analyze
specific statistical models and certify that the RE
condition (\ref{EqnRSCConstants}) and deviation conditions are
satisfied by a random pair $(\GamHat, \widehat{\gamma
})$ with high probability.
We note that for the standard Lasso choice $(\GamLas,
\widehat{\gamma}_{\mathrm{Las}})$ of
this matrix--vector pair, bounds of the form (\ref{EqnGenBound}) for
sub-Gaussian noise are well known from past work
(e.g.,~\cite{BicEtal08,Huang06,MeiYu09,Neg09}). The novelty of
Theorem~\ref{ThmStat} is in allowing for general pairs of such
surrogates, which---as shown by the examples discussed earlier---can
lead to nonconvexity in the underlying $M$-estimator. Moreover, some
interesting differences arise due to the term $\varphi(\Qprob,
\sigeps)$, %which changes depending on the nature of the model
%(missing, noisy, and/or dependent), as will be clarified in the
%sequel, proving that the conditions of Theorem~\ref{ThmStat} are
%satisfied with high probability, for noisy/missing data requires some
%nontrivial analysis, involving both concentration inequalities and
%random matrix theory.
which changes depending on the nature of the model (missing, noisy,
and/or dependent). As will be clarified in the sequel. Proving that the
conditions of Theorem~\ref{ThmStat} are satisfied with high probability for
noisy/missing data requires some nontrivial analysis involving both
concentration inequalities and random matrix theory.

Note that in the presence of nonconvexity, it is possible in
principle for the optimization problems (\ref{EqnGeneral})
and (\ref{EqnGeneralSlight}) to have \textit{many} global optima that
are separated by large distances. Interestingly, Theorem~\ref{ThmStat}
guarantees that this unpleasant feature does not arise under the
stated conditions: given any two global optima $\betahat$ and
$\betatil$ of the program (\ref{EqnGeneral}), Theorem~\ref{ThmStat}
combined with the triangle inequality guarantees that
\[
\|\betahat- \betatil\|_2 \leq\bigl\|\betahat- \betastar\bigr\|_2 +
\bigl\|\betatil- \betastar\bigr\|_2 \leq2 \plaincon_0
\frac{\varphi
(\Qprob,
\sigeps)}{\allow} \sqrt{\frac{\kdim\log\pdim}{\numobs}}
\]
[and similarly for the
program (\ref{EqnGeneralSlight})]. Consequently, under any scaling
such that $\frac{\kdim\log\pdim}{\numobs} = o(1)$, the set of all
global optima must lie
within an $\ell_2$-ball whose radius shrinks to zero.

%f1 #&#
\begin{figure}
\begin{tabular}{@{}cc@{}}

\includegraphics{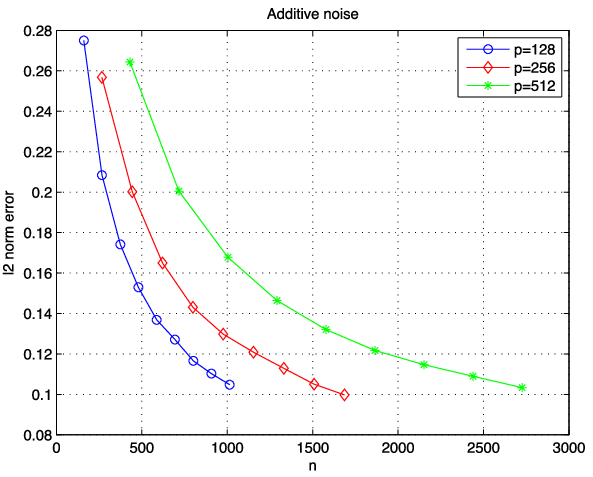}
 & \includegraphics{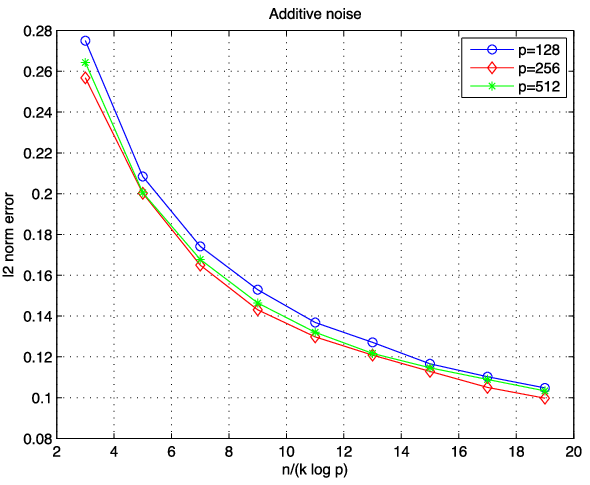}\\
(a) & (b)
\end{tabular}
\caption{Plots of the error $\|\betahat- \betastar\|_2$ after running
projected gradient descent on the nonconvex objective, with
sparsity $k \approx\sqrt{p}$. Plot \textup{(a)} is an error plot for
i.i.d. data with additive noise, and plot \textup{(b)} shows $\ell_2$-error
versus the rescaled sample size $\frac{n}{k \log p}$. As predicted
by Theorem \protect\ref{ThmStat}, the curves align for different
values of
$p$ in the rescaled plot.}
\label{FigRescaling}
\end{figure}

In addition, it is worth observing that Theorem~\ref{ThmStat} makes a
specific prediction for the scaling behavior of the $\ell_2$-error
$\|\betahat- \betastar\|_2$. In order to study this scaling
prediction, we performed simulations under the additive noise model
described in Example~\ref{ExaNoisy}, using the parameter setting
$\Sigma_x = I$ and $\Sigma_w = \sigma_w^2 I$ with $\sigma_w = 0.2$.
Panel (a) of Figure~\ref{FigRescaling} provides
plots\setcounter{footnote}{2}\footnote{Corollary~\ref{CorAddNoise}, to be stated shortly,
guarantees that the conditions of Theorem~\ref{ThmStat} are
satisfied with high probability for the additive noise model. In
addition, Theorem~\ref{ThmOpt} to follow provides an efficient
method of obtaining an accurate approximation of the global
optimum.} of the error $\|\betahat- \betastar\|_2$ versus the
sample size~$\numobs$, for problem dimensions $\pdim\in\{128,
256, 512 \}$. Note that for all three choices of dimensions, the
error decreases to zero as the sample size $\numobs$ increases,
showing consistency of the method. The curves also shift to the right
as the dimension $p$ increases, reflecting the natural intuition that
larger problems are harder in a certain sense. Theorem~\ref{ThmStat}
makes a~specific prediction about this scaling behavior: in
particular, if we plot the $\ell_2$-error versus the rescaled sample
size $\numobs/(\kdim\log\pdim)$, the curves should roughly align for
different values of $p$. Panel (b) shows the same data re-plotted on
these rescaled axes, thus verifying the predicted ``stacking
behavior.''

Finally, as noted by a reviewer, the constraint $R = \|\betastar\|_1$
in the program~(\ref{EqnGeneral}) is rather restrictive, since
$\betastar$ is unknown. Theorem~\ref{ThmStat} merely establishes a
heuristic for the scaling expected for this optimal radius. In this
regard, the Lagrangian estimator (\ref{EqnGeneralSlight}) is more
appealing, since it only requires choosing $b_0$ to be larger than
$\|\betastar\|_2$, and the conditions on the regularizer $\regpar$ are
the standard ones from past work on the Lasso.

%s3.1.2 #&#
\subsubsection{Optimization error}

Although Theorem~\ref{ThmStat} provides guarantees that hold uniformly
for any global minimizer, it does not provide guidance on how to
approximate such a global minimizer using a polynomial-time algorithm.
Indeed, for nonconvex programs in general, gradient-type methods may
become trapped in local minima, and it is impossible to guarantee that
all such local minima are close to a global optimum. Nonetheless, we
are able to show that for the family of programs (\ref{EqnGeneral}),
under reasonable conditions on $\GamHat$ satisfied in various
settings, simple gradient methods will converge geometrically fast to
a very good approximation of any global optimum. The following
theorem supposes that we apply the projected gradient
updates~(\ref{EqnProjGrad}) to the constrained
program~(\ref{EqnGeneral}), or the composite updates (\ref{EqnMirror})
to the Lagrangian program (\ref{EqnGeneralSlight}), with stepsize
$\eta= 2 \alup$. In both cases, we assume that $\numobs
\succsim\kdim\log\pdim$, as is required for statistical consistency
in Theorem~\ref{ThmStat}.

%th2 #&#
\begin{theos}[(Optimization error)]
\label{ThmOpt}
Under the conditions of Theorem~\ref{ThmStat}:
\begin{longlist}[(a)]
\item[(a)] For any global optimum $\betahat$ of the constrained
program (\ref{EqnGeneral}), there are universal positive constants
$(\plaincon_1, \plaincon_2)$ and a contraction coefficient $\CONTRAC
\in(0,1)$, independent of $(n,p,k)$, such that the gradient descent
iterates (\ref{EqnProjGrad}) satisfy the bounds
%
%e3.5 #&#
%e3.6 #&#
\begin{eqnarray}
\label{EqnOpt} \bigl\|\beta^t - \betahat\bigr\|_2^2
& \leq &\CONTRAC^t \bigl\|\beta^0 - \betahat\bigr\|_2^2
+ \plaincon_1 \frac{\log\pdim}{\numobs} \bigl\|\betahat- \betastar
\bigr\|_1^2 + \plaincon_2 \bigl\|\betahat- \betastar
\bigr\|_2^2,
\\
\label{EqnOptEll1} \bigl\|\beta^t - \betahat\bigr\|_1
& \leq & 2 \sqrt{\kdim} \bigl\|\beta^t - \betahat\bigr\|_2 + 2 \sqrt{
\kdim} \bigl\|\betahat- \betastar\bigr\|_2 + 2 \bigl\|\betahat- \betastar
\bigr\|_1
\end{eqnarray}
for all $t \ge0$.
\item[(b)] Letting $\phi$ denote the objective function of Lagrangian
program (\ref{EqnGeneralSlight}) with global optimum $\betahat$, and
applying composite gradient updates (\ref{EqnMirror}), there are
universal positive constants $(\plaincon_1, \plaincon_2)$ and a
contraction coefficient $\CONTRAC\in(0,1)$, independent of
$(\numobs, \pdim, \kdim)$, such that
%
%e3.7 #&#
\begin{equation}
\label{EqnOptCon} \bigl\|\beta^t - \betahat\bigr\|_2^2
\leq\underbrace{\plaincon_1 \bigl\| \betahat- \betastar
\bigr\|_2^2}_{\delta^2} \qquad\mbox{for all iterates $t \geq
T$,}
\end{equation}
where $T \defn c_2\log\frac{(\phi(\beta^0) -
\phi(\betahat))}{\delta^2}  /\log(1/\CONTRAC)$.
\end{longlist}
\end{theos}

\textit{Remarks.\quad}
As with Theorem~\ref{ThmStat}, these claims are deterministic in
nature. Probabilistic conditions will enter into the corollaries, which
involve proving that the surrogate matrices $\GamHat$ used for noisy,
missing and/or dependent data satisfy the lower- and upper-RE
conditions with high probability. The proof of Theorem~\ref{ThmOpt}
itself is based on an extension of a result due to Agarwal et
al.~\cite{AgaEtal11} on the convergence of projected gradient descent
and composite gradient descent in high dimensions. Their result, as
originally stated, imposed convexity of the loss function, but the
proof can be modified so as to apply to the nonconvex loss functions
of interest here. As noted following Theorem~\ref{ThmStat},
all global minimizers of the nonconvex program (\ref{EqnGeneral})
lie within a small ball. In addition, Theorem~\ref{ThmOpt}
guarantees that the local minimizers also lie within a ball of the
same magnitude. Note that in order to show that
Theorem~\ref{ThmOpt} can be applied to the specific statistical models
of interest in this paper, a considerable amount of technical analysis
remains in order to establish that its conditions hold with high
probability.

In order to understand the significance of the bounds (\ref{EqnOpt})
and (\ref{EqnOptCon}), note that they provide upper bounds for the
$\ell_2$-distance between the iterate~$\beta^t$ at time $t$, which is
easily computed in polynomial-time, and any global optimum $\betahat$
of the program (\ref{EqnGeneral}) or (\ref{EqnGeneralSlight}), which
may be difficult to compute. Focusing on bound (\ref{EqnOpt}),
since $\gamma\in(0,1)$, the first term in the bound vanishes as $t$
increases. The remaining terms involve the statistical errors
$\|\betahat- \betastar\|_q$, for $q = 1, 2$, which are controlled in
Theorem~\ref{ThmStat}. It can be verified that the two terms
involving the statistical error on the right-hand side are bounded as
$\order(\frac{\kdim\log\pdim}{\numobs})$, so Theorem~\ref{ThmOpt}\vadjust{\goodbreak}
guarantees that projected gradient descent produce an output that is
essentially as good---in terms of statistical error---as any global
optimum of the program (\ref{EqnGeneral}). Bound (\ref{EqnOptCon})
provides a~similar guarantee for composite gradient descent applied to
the Lagrangian version.
%Similarly,
%inequality \eqref{EqnOptEll1}---when combined with the statistical
%bounds from Theorem~\ref{ThmStat}---yields that the optimization error
%in $\ell_1$-norm is bounded as $\order\big(k\sqrt{\frac{\log
% p}{n}}\big)$, which is of the same order as the
%$\ell_1$-statistical error.

%f2 #&#
\begin{figure}
\begin{tabular}{@{}cc@{}}

\includegraphics{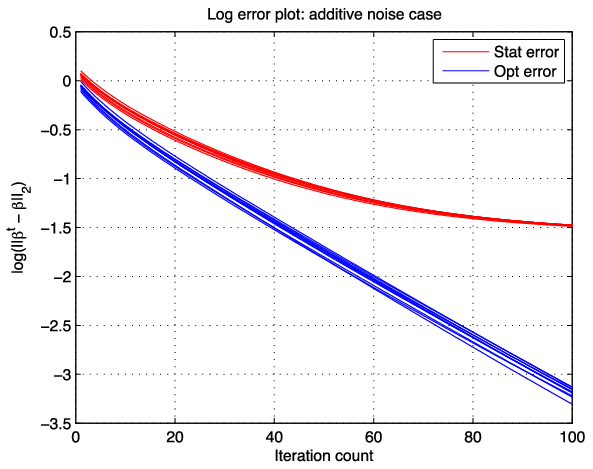}
 & \includegraphics{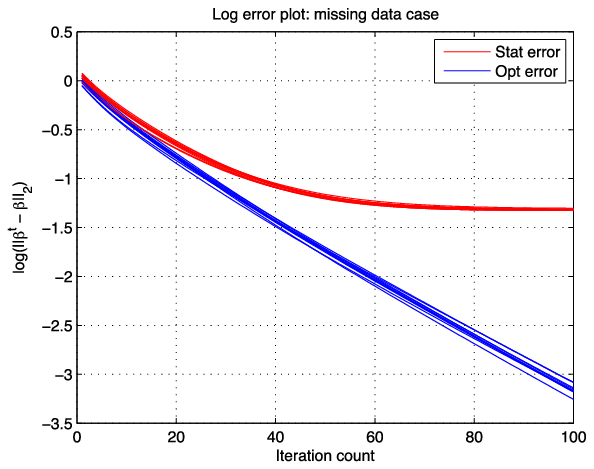}\\
(a) & (b)\vspace*{-3pt}
\end{tabular}
\caption{Plots of the optimization error $\log(\|\beta^t -
\betahat\|_2)$ and statistical error $\log(\|\beta^t -
\betastar\|_2)$ versus iteration number $t$, generated by running
projected gradient descent on the nonconvex objective. Each plot
shows the solution path for the same problem instance, using $10$
different starting points. As predicted by Theorem \protect\ref
{ThmOpt}, the
optimization error decreases geometrically.}
\label{FigOptErr}\vspace*{-3pt}
\end{figure}

Experimentally, we have found that the predictions of
Theorem~\ref{ThmOpt} are borne out in simulations.
Figure~\ref{FigOptErr} shows the results of applying the projected
gradient descent method to solve the optimization
problem (\ref{EqnGeneral}) in the case of additive noise [panel (a)],
and missing data [panel (b)]. In each case, we generated a random
problem instance, and then applied the projected gradient descent
method to compute an estimate $\betahat$. We then reapplied the
projected gradient method to the same problem instance $10$ times,
each time with a random starting point, and measured the error
$\|\beta^t - \widehat{\beta}\|_2$ between the iterates and the first
estimate (optimization error), and the error $\|\beta^t -
\betastar\|_2$ between the iterates and the truth (statistical error).
Within each panel, the blue traces show the optimization error over
$10$ trials, and the red traces show the statistical error. On the
logarithmic scale given, a~geometric rate of convergence corresponds
to a straight line. As predicted by Theorem~\ref{ThmOpt}, regardless
of the starting point, the iterates $\{\beta^t \}$ exhibit geometric
convergence to the same fixed point.\footnote{To be precise,
Theorem~\ref{ThmOpt} states that the iterates will converge
geometrically to a small neighborhood of all the global optima.}
The statistical error contracts geometrically up to a certain point,
then flattens out.

%%%%%%%%%%%%%%%%%%%%%%%%%%%%%%%%%%%%%%%%%%%%%%%%%%%%%%%%%%%%%%%%%%%%%%%%%%%%%%%

%s3.2 #&#
\subsection{Some consequences}
\label{SecConseq}

As discussed previously, both Theorems~\ref{ThmStat} and~\ref{ThmOpt}
are deterministic results. Applying them to specific statistical
models requires some\vadjust{\goodbreak} additional work in order to establish that the
stated conditions are met. We now turn to the statements of some
consequences of these theorems for different cases of noisy, missing
and dependent data. In all the corollaries below, the claims hold
with probability greater than $1 - \plaincon_1 \exp(-\plaincon_2
\log
\pdim)$, where $(\plaincon_1, \plaincon_2)$ are universal positive
constants, independent of all other problem parameters. Note that in
all corollaries, the triplet $(\numobs, \pdim, \kdim)$ is assumed to
satisfy scaling of the form $\numobs\succsim\kdim\log\pdim$, as is
necessary for $\ell_2$-consistent estimation of $\kdim$-sparse vectors
in $\pdim$ dimensions.\looseness=1
%
%de3 #&#
\begin{defns}
\label{DefnSubGaussMat}
We say that a random matrix $X \in\real^{\numobs\times\pdim}$ is
sub-Gaussian with parameters $(\Sigma,\sigma^2)$ if:
\begin{longlist}[(a)]
\item[(a)] each row $x_i^T \in\real^\pdim$ is sampled independently
from a zero-mean distribution with covariance $\Sigma$, and
\item[(b)] for any unit vector $u \in\real^\pdim$, the random
variable $u^T x_i$ is sub-Gaussian with parameter at most $\sigma$.
\end{longlist}
\end{defns}

For instance, if we form a random matrix by drawing each row
independently from the distribution $N(0, \Sigma)$, then the resulting
matrix $X \in\real^{\numobs\times\pdim}$ is a sub-Gaussian matrix
with parameters $(\Sigma, \Vvert\Sigma\Vvert_{\mathrm{op}})$.

%s3.2.1 #&#
\subsubsection{Bounds for additive noise: i.i.d. case}

We begin with the case of i.i.d. samples with additive noise, as
described in Example~\ref{ExaNoisy}.
%
%co1 #&#
\begin{cors}
\label{CorAddNoise}
Suppose that we observe $Z = X + W$, where the random matrices $X,
W\in\real^{n \times p}$ are sub-Gaussian with parameters $(\Sigma_x,
\sigma_x^2)$, and let~$\varepsilon$ be an i.i.d. sub-Gaussian
vector with parameter $\sigeps^2$. Let $\sigma_z^2 = \sigma_x^2 +
\sigma_w^2$. Then under the scaling $n \succsim
\max \{\frac{\sigma_z^4}{\lammin^2(\Sigma_x)}, 1 \} k \log p$,
for the $M$-estimator based on the surrogates $(\GamAdd,
\widehat{\gamma}_{\mathrm{add}})$,
the results of Theorems~\ref{ThmStat} and~\ref{ThmOpt} hold with
parameters $\allow= \frac{1}{2} \lambda_{\min}(\Sigma_x)$ and
$\varphi(\Qprob, \sigeps) = \plaincon_0 \sigma_z(\sigma_w +
\sigeps)
\|\betastar\|_2$, with probability at least $1 - \plaincon_1
\exp(-\plaincon_2 \log\pdim)$. %\pldel{and $(\Sigma_w,
%Then for the $M$-estimator based on the surrogates $(\GamAdd,
%the results of Theorems~\ref{ThmStat} and~\ref{ThmOpt}
%hold with parameters
%%
%
%%
%}
%
\end{cors}

\textit{Remarks.\quad}
(a) Consequently, the $\ell_2$-error of any optimal solution
$\betahat$ satisfies the bound
\[
\bigl\|\betahat- \betastar\bigr\|_2 \precsim\frac{\sigma_z(\sigma_w +
\sigeps)}{\lammin(\Sigma_x)} \bigl\|\betastar
\bigr\|_2 \sqrt{\frac{k \log
p}{n}}
\]
%
%%
%}{\lammin(\Sigma_x)}
% \sqrt{\frac{\kdim\log\pdim}{\numobs}}
%%
%}
with high probability. The prefactor in this bound has a natural
interpretation as an inverse signal-to-noise ratio; for instance, when
$X$ and $W$ are zero-mean Gaussian matrices with row covariances
$\Sigma_x = \sigma_x^2 I$ and $\Sigma_w = \sigma_w^2 I$, respectively,
we have $\lammin(\Sigma_x) = \sigma_x^2$, so
\[
\frac{(\sigma_w + \sigeps) \sqrt{\sigma_x^2 +
\sigma_w^2}}{\lammin(\Sigma_x)} = \frac{\sigma_w +
\sigeps}{\sigma_x} \sqrt{1 + \frac{\sigma_w^2}{\sigma_x^2}}.
\]
This quantity grows with the ratios $\sigma_w/\sigma_x$ and
$\sigeps/\sigma_x$, which measure the SNR of the observed covariates
and predictors, respectively. Note that when $\sigma_w = 0$,
corresponding to the case of uncorrupted covariates, the bound on
$\ell_2$-error agrees with known results. %\pldel{
%%
%{1 +
%%
%This quantity grows with the ratio $\frac{\sigma_w^2}{\sigma_x^2}$,
%which measures noisiness in the covariates, whereas the quantity
%$\sigeps$ is the standard deviation of the usual observation noise.}
See Section~\ref{SecSims} for simulations and further discussions of
the consequences of Corollary~\ref{CorAddNoise}.

(b) We may also compare the results in (a) with bounds from
past work on high-dimensional sparse regression with noisy
covariates~\cite{RosTsy11}. In this work, Rosenbaum and
Tsybakov derive similar concentration bounds on sub-Gaussian
matrices. The tolerance parameters are all
$\order (\sqrt{\frac{\log p}{n}} )$, with prefactors
depending on the sub-Gaussian parameters of the matrices. In
particular, in their notation,
\[
\nu\asymp\bigl(\sigma_x \sigma_w +
\sigma_w \sigeps+ \sigma_w^2\bigr) \sqrt{
\frac{\log p}{n}} \bigl\|\betastar\bigr\|_1 ,
\]
leading to the bound (cf. Theorem 2 of Rosenbaum and
Tsybakov~\cite{RosTsy11})
\[
\bigl\|\betahat- \betastar\bigr\|_2 \precsim \frac{\nu\sqrt{k}}{\lammin(\Sigma_x)} \asymp
\frac{\sigma^2}{\lammin(\Sigma_x)}  \sqrt {\frac{
k \log p}{n}}\bigl\|\betastar\bigr\|_1.
\]
%
%which is larger than our bound by a factor of $\sqrt{k}$. %\pldel{In
%this work, Rosenbaum and Tsybakov require an upper bound on the
%matrix $W$ of the form $\|W\|_{\max} < \delta$ to hold with high
%probability, for some constant $\delta> 0$. Since $W$ is a matrix
%of i.i.d. Gaussians with variance $\sigma_w^2$, we have
%$\|W\|_{\max} \asymp\sigma_w \sqrt{\log(pn)} \approx\sigma_w
%Tsybakov~\cite{RosTsy11} (and assuming $\sigeps= 0$) then gives
%%
%%
%compared to the sharper bound
%%
%%
%from our theorem.}
%
\vspace*{8pt}

\textit{Extensions to unknown noise covariance.\quad} Situations may
arise where the noise covariance $\Sigma_w$ is unknown, and must be
estimated from the data. One simple method is to assume that $\Sigma_w$
is estimated from independent observations of the noise. In this case,
suppose we independently observe\vspace*{1pt} a~matrix $W_0 \in\real^{n
\times p}$ with $n$ i.i.d. vectors of noise. Then we use $\SigHat_w =
\frac{1}{n} W_0^TW_0$ as our\vspace*{1pt} estimate of $\Sigma_w$. A more
sophisticated variant of this method (cf. Chapter 4 of Carroll et al.
\cite{CarEtal95}) assumes that we observe $k_i$ replicate measurements
$Z_{i1}, \ldots, Z_{ik}$ for each $x_i$ and form the estimator
%
%e3.8 #&#
\begin{equation}
\label{EqnSigHat} \SigHat_w = \frac{\sum_{i=1}^n \sum_{j=1}^{k_i} (Z_{ij} -
\Zbar_{i\cdot}) (Z_{ij} - \Zbar_{i\cdot})^T}{\sum_{i=1}^n (k_i -
1)}.
\end{equation}
Based on the estimator $\SigHat_w$, we form the pair $(\GamTil,
\widetilde{\gamma})$ such that $\widetilde
{\gamma}= \frac{1}{n} Z^T y$ and $\GamTil=
\frac{Z^TZ}{n} - \SigHat_w$. In the proofs of Section~\ref{SecProofs},
we will analyze the case where $\SigHat_w = \frac{1}{n} W_0^TW_0$ and
show that the result of Corollary~\ref{CorAddNoise} still holds when
$\Sigma_w$ must be estimated from the data. Note that the estimator in
equation (\ref{EqnSigHat}) will also yield the same result, but the
analysis is more complicated.

%s3.2.2 #&#
\subsubsection{Bounds for missing data: i.i.d. case}

Next, we turn to the case of i.i.d. samples with missing data, as
discussed in Example~\ref{ExaMul}. For a missing data parameter vector
$\rhobf$, we define $\rho_{\max} \defn\max_j \rho_j$, and assume
$\rho_{\max} < 1$.
%
%co2 #&#
\begin{cors}
\label{CorMAR}
Let $X \in\real^{n \times p}$ be sub-Gaussian with parameters
$(\Sigma_x, \sigma_x^2)$, and $Z$ the missing data matrix with
parameter $\rhobf$. Let $\varepsilon$ be an\vadjust{\goodbreak} i.i.d. sub-Gaussian
vector with parameter $\sigeps^2$. If $n \succsim
\max (\frac{1}{(1 - \rho_{\max})^4}
\frac{\sigma_x^4}{\lambda_{\min}^2(\Sigma_x)}, 1 ) k \log p$, then
Theorems~\ref{ThmStat} and~\ref{ThmOpt} hold with probability at least
$1 - c_1 \exp(-c_2 \log p)$ for $\allow= \frac{1}{2}
\lambda_{\min}(\Sigma_x)$ and $\varphi(\Qprob, \sigma_\varepsilon)
= c_0
\frac{\sigma_x}{1-\rho_{\max}}  (\sigeps+
\frac{\sigma_x}{1-\rho_{\max}} ) \|\betastar\|_2$.
\end{cors}

\textit{Remarks.\quad}
Suppose $X$ is a Gaussian random matrix and $\rho_j = \rho$
for all $j$. In this case, the ratio
$\frac{\sigma_x^2}{\lambda_{\min}(\Sigma_x)} =
\frac{\lambda_{\max}(\Sigma_x)}{\lambda_{\min}(\Sigma_x)} =
\kappa(\Sigma_x)$ is the condition number of $\Sigma_x$. Then
\[
\frac{\varphi(\Qprob, \sigma_\varepsilon)}{\alpha} \asymp \biggl(\frac{1}{\lammin(\Sigma_x)} \frac{\sigma_x\sigeps}{1-\rho
} +
\frac{\kappa(\Sigma_x)}{(1-\rho)^2} \biggr) \bigl\|\betastar\bigr\|_2,
\]
a quantity that depends on both the conditioning of $\Sigma_x$, and
the fraction $\rho\in[0,1)$ of missing data. We will consider the
results of Corollary~\ref{CorMAR} applied to this example in the
simulations of Section~\ref{SecSims}.

%(b) Again, it is instructive to compare our result with the
%past work of Rosenbaum and Tsybakov~\cite{RosTsy11}. In the
%case of missing data, their estimator is not applicable the case
%of random design $X$, since it models missing data as the additive
%model $X+W$, where $W$ is an independent noise term. An earlier
%paper by the same authors~\cite{RosTsy10} is applicable to missing
%data in random designs, but requires that $\|X\|_{\max}$ is
%bounded by some universal constant with high probability. Using
%standard results on the maximum of i.i.d. Gaussians, their
%estimator is guaranteed to return an estimate such that $\|\betahat
%- \betastar\|_2 \precsim\frac{\sigma_x}{\sigma_x} \sqrt{k \log
%larger than our bound by a factor of $\sqrt{\numobs}$. %\pldel{A more
%%refined argument, stated without proof in their paper, leads to
%%the guarantee \mbox{$\|\betahat- \betastar\|_2 \precsim\sqrt{k}
%%\big(\rhoprob+ \sqrt{\frac{\log p}{n}} \big)$.} Although this
%%bound improves the initial bound \eqref{EqnRosTsyInitial}, it is
%%still weaker than our guarantees; in particular, it does not
%%decrease to zero as the sample size $\numobs$ tends to infinity,
%%and it is much larger whenever $\numobs\gg\log\pdim$, the
%%regime of interest.}
%%\end{longlist}
%
\vspace*{8pt}

\textit{Extensions to unknown $\rho$.\quad}
As in the additive noise case, we may wish to consider the case when
the missing data parameters~$\rhobf$ are not observed and must be
estimated from the data. For each $j = 1, 2, \ldots, \pdim$, we
estimate~$\rhoprob_j$ using $\rhoprobhat_j$, the empirical average of
the number of observed entries per column. Let $\rhohatbf\in\real^p$
denote the resulting estimator of $\rhobf$. Naturally, we use the pair
of estimators $(\GamTil, \widetilde{\gamma})$ defined by
%
%e3.9 #&#
\begin{equation}
\label{EqnDefnGamUnknown} \GamTil= \frac{Z^TZ}{n} \hackdiv\Mtil
\quad\mbox{and}\quad
\widetilde{\gamma}= \frac{1}{\numobs} Z^T y \hackdiv(\mathbf{1} -
\rhohatbf),
\end{equation}
where
\[
\Mtil_{ij} = %
\cases{ (1-\rhoprobhat_i) (1-
\rhoprobhat_j), &\quad if $i \neq j$,
\cr
1-
\rhoprobhat_i, &\quad if $i = j$.} %
\]
We will show in Section~\ref{SecProofs} that Corollary~\ref{CorMAR}
holds when $\rhobf$ is estimated by $\rhohatbf$.

%%%%%%%%%%%%%%%%%%%%%%%%%%%%%%%%%%%%%%%%%%%%%%%%%%%%%%%%%%%%%%%%%%%%%%%%

%s3.2.3 #&#
\subsubsection{Bounds for dependent data}

Turning to the case of dependent data, we consider the setting where
the rows of $X$ are drawn from a stationary vector autoregressive
(VAR) process according to
%
%e3.10 #&#
\begin{equation}
\label{EqnVAR} x_{i+1} = A x_i + v_i
\qquad\mbox{for $i = 1, 2, \ldots, \numobs -1$},
\end{equation}
where $v_i \in\real^\pdim$ is a zero-mean noise vector with
covariance matrix $\CovV$, and $A \in\real^{\pdim\times\pdim}$ is a
driving matrix with spectral norm $\Vvert A \Vvert_{{2}} < 1$. We
assume the
rows of $X$ are drawn from a Gaussian distribution with covariance
$\Sigma_x$, such that $\Sigma_x = A\Sigma_x A^T + \Sigma_v$. Hence,
the rows of $X$ are identically distributed but not independent, with
the choice $A = 0$ giving rise to the
i.i.d. scenario. Corollaries~\ref{CorARaddNoise}
and~\ref{CorARmissing} correspond to the case of additive noise and
missing data for a Gaussian VAR process.

%
%co3 #&#
\begin{cors}
\label{CorARaddNoise}
Suppose the rows of $X$ are drawn according to a \mbox{Gaussian} VAR process
with driving matrix $A$. Suppose the additive noise matrix~$W$ is\vadjust{\goodbreak}
i.i.d. with Gaussian rows, and let $\varepsilon$ be an i.i.d.
sub-Gaussian vector with parameter~$\sigeps^2$. If $n
\succsim\max (\frac{\zeta^4}{\lambda_{\min}^2(\Sigma_x)},
1 )
k \log p$, with $\zeta^2 = \Vvert\Sigma_w \Vvert_{\mathrm{op}} +
\frac{2\Vvert\Sigma_x \Vvert_{\mathrm{op}}}{1 - \Vvert A \Vvert_{\mathrm{op}}}$, then
Theorems~\ref{ThmStat} and~\ref{ThmOpt} hold with probability at least
$1 - c_1 \exp(-c_2 \log p)$ for $\allow= \frac{1}{2}
\lambda_{\min}(\Sigma_x)$ and $\varphi(\Qprob, \sigma_\varepsilon) =
c_0(\sigeps\zeta+ \zeta^2) \|\betastar\|_2$.
\end{cors}

%co4 #&#
\begin{cors}
\label{CorARmissing}
Suppose the rows of $X$ are drawn according to a Gaussian VAR process
with driving matrix $A$, and $Z$ is the observed matrix subject to
missing data, with parameter $\rhobf$. Let $\varepsilon$ be an
i.i.d. sub-Gaussian vector with parameter~$\sigeps^2$. If $n
\succsim\max (\frac{\zeta^{\prime
4}}{\lambda_{\min}^2(\Sigma_x)}, 1 ) k \log p$, with
$\zeta^{\prime2} =
\frac{1}{(1-\rhomax)^2}\frac{2\Vvert\Sigma_x \Vvert_{\mathrm
{op}}}{1 - \Vvert A \Vvert_{\mathrm{op}}}$,
then Theorems~\ref{ThmStat} and~\ref{ThmOpt} hold with probability at
least $1 - c_1 \exp(-c_2 \log p)$ for $\allow= \frac{1}{2}
\lambda_{\min}(\Sigma_x)$ and $\varphi(\Qprob, \sigma_\varepsilon) =
c_0(\sigeps\zeta^\prime+ \zeta^{\prime2}) \|\betastar\|_2$.
\end{cors}

\begin{rem*}
Note that the scaling and the form of $\varphi$ in Corollaries
\ref{CorMAR}--\ref{CorARmissing}
are very similar, except with different effective\vspace*{1pt} variances
$\sigma^2 = \frac{\sigma_x^2}{(1-\rhomax)^2}$, $\zeta^2$ or
$\zeta^{\prime2}$, depending on the type of corruption in the
data. As we will see in Section~\ref{SecProofs}, the proofs involve
verifying the deviation conditions (\ref{EqnDev}) using similar
techniques. On the other hand, the proof of Corollary~\ref{CorAddNoise} proceeds via
deviation condition (\ref{EqnDevFull}), which produces a tighter
bound.

Note that we may extend the cases of dependent data to situations
when~$\Sigma_w$ and $\rhobf$ are unknown and must be estimated from
the data. The proofs of these extensions are identical to the i.i.d
case, so we will omit them.
\end{rem*}

%%%%%%%%%%%%%%%%%%%%%%%%%%%%%%%%%%%%%%%%%%%%%%%%%%%%%%%%%%%%%%%%%%%%%%%%%%%%%%%

%s3.3 #&#
\subsection{Application to graphical model inverse covariance estimation}

The problem of inverse covariance estimation for a Gaussian graphical
model is also related to the Lasso. Meinshausen and
B\"{u}hlmann~\cite{MeiBuh06} prescribed a~way to recover the support
of the precision matrix $\Theta$ when each column of $\Theta$ is
$k$-sparse, via linear regression and the Lasso. More recently,
Yuan~\cite{Yua10} proposed a method for estimating $\Theta$ using the
Dantzig selector, and obtained error bounds on $\Vvert\Thetahat-
\Theta\Vvert_{{1}}$ when the columns of $\Theta$ are bounded in
$\ell_1$. Both of these results assume that $X$ is fully-observed and has
i.i.d. rows.

Suppose we are given a matrix $X \in\real^{n \times p}$ of samples
from a multivariate Gaussian distribution, where each row is
distributed according to $N(0, \Sigma)$. We assume the rows of $X$ are
either i.i.d. or sampled from a Gaussian VAR process. Based on the
modified Lasso of the previous section, we devise a~method to estimate~$\Theta$
based on a corrupted observation matrix $Z$, when~$\Theta$ is
sparse. Our method bears similarity to the method of Yuan~\cite{Yua10},
but is valid in the case of corrupted data, and does not require an
$\ell_1$ column bound. Let~$X^j$ denote the $j$th column
of~$X$, and let $X^{-j}$ denote the matrix~$X$ with $j$th
column removed. By standard results on Gaussian graphical models,
there exists a vector $\theta^j \in\real^{p-1}$ such that
%
%e3.11 #&#
\begin{equation}
\label{EqnGaussReg} X^j = X^{-j} \theta^j +
\varepsilon^j,
\end{equation}
where $\varepsilon^j$ is a vector of i.i.d. Gaussians and $\varepsilon^j
\condind X^{-j}$ for each $j$. If we define $a_j \defn
-(\Sigma_{jj} - \Sigma_{j, -j} \theta^j)^{-1}$, we can verify
that\vadjust{\goodbreak}
$\Theta_{j,-j} = a_j \theta^j$. Our algorithm, described below, forms
estimates $\widehat{\theta}{}^j$ and $\ahat_j$ for each
$j$, then combines the
estimates to obtain an estimate $\Thetahat_{j, -j} = \ahat_j
\widehat{\theta}{}^j$.

In the additive noise case, we observe the matrix $Z = X+W$. From the
equations (\ref{EqnGaussReg}), we obtain $Z^j = X^{-j}\theta^j +
(\varepsilon^j + W^j)$. Note that $\delta^j = \varepsilon^j + W^j$ is a
vector of i.i.d. Gaussians, and since $X \condind W$, we have
$\delta^j \condind X^{-j}$. Hence, our results on covariates with
additive noise allow us to recover $\theta^j$ from $Z$. We can verify
that this reduces to solving the\vspace*{2pt} program (\ref{EqnGeneral})
or (\ref{EqnGeneralSlight}) with the pair $(\GamHat^{(j)},
\widehat{\gamma}{}^{(j)}) = (\SigHat_{-j,-j}, \frac
{1}{n} Z^{-j T} Z^j)$,
where $\SigHat= \frac{1}{n} Z^TZ - \Sigma_w$.

When $Z$ is a missing-data version of $X$, we similarly estimate the
vectors~$\theta^j$ via equation (\ref{EqnGaussReg}), using our results
on the Lasso with missing covariates. Here, both covariates and
responses are subject to missing data, but this makes no difference in
our theoretical results. For each $j$, we use the pair
\[
\bigl(\GamHat^{(j)}, \widehat{\gamma}{}^{(j)}\bigr) = \biggl(
\SigHat_{-j,-j}, \frac{1}{n} Z^{-j T}Z^j \hackdiv
\bigl(\mathbf{1}-\rhobf^{-j}\bigr) (1-\rho_j) \biggr),
\]
where $\SigHat= \frac{1}{n} Z^TZ \hackdiv M$, and $M$ is defined as in
Example~\ref{ExaMul}.

To obtain the estimate $\Thetahat$, we therefore propose the following
procedure, based on the estimators $\{(\GamHat^{(j)},
\widehat{\gamma}{}^{(j)})\}_{j=1}^p$ and $\SigHat$.
%
%al3.1 #&#
\begin{alg}
\label{AlgMatrix}
(1) Perform $\pdim$ linear regressions of the variables $Z^j$
upon the remaining variables $Z^{-j}$, using\vspace*{1pt} the
program (\ref{EqnGeneral}) or (\ref{EqnGeneralSlight}) with the
estimators $(\GamHat^{(j)}, \widehat{\gamma}{}^{(j)})$,
to obtain estimates
$\widehat{\theta}{}^j$ of $\theta^j$.

(2) Estimate\vspace*{1pt} the scalars $a_j$ using the quantity $\ahat_j
\defn-(\SigHat_{jj} - \SigHat_{j,-j}\widehat{\theta
}{}^j)^{-1}$, based on
the estimator $\SigHat$. Form $\Thetatil$ with $\Thetatil_{j,-j}
= \ahat_j \widehat{\theta}{}^j$ and $\Thetatil_{jj} =
-\ahat_j$.\vspace*{2pt}

(3) Set $\Thetahat= \arg\min_{\Theta\in S^p}
\Vvert\Theta- \Thetatil\Vvert_{{1}}$, where $S^\pdim$ is the set of
symmetric matrices.
\end{alg}

Note that the minimization in step (3) is a linear program,
so is easily solved with standard methods. We have the following
corollary about $\Thetahat$:

%f3 #&#
\begin{figure}[b]
\begin{tabular}{@{}cc@{}}

\includegraphics{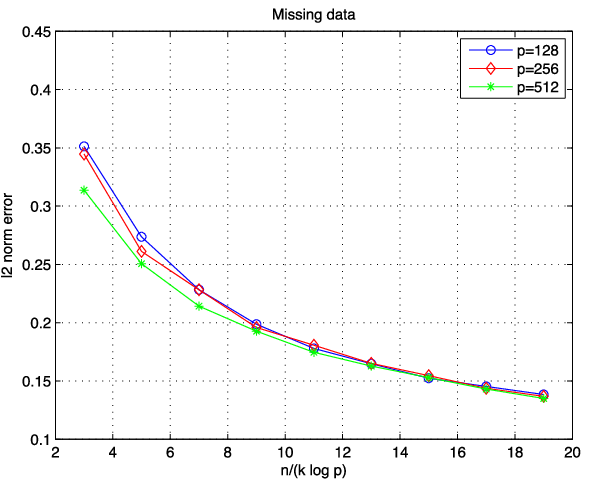}
 & \includegraphics{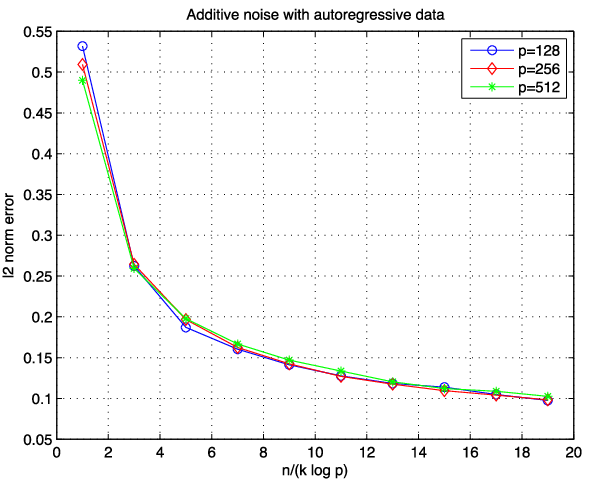}\\
(a) & (b)
\end{tabular}
\caption{Plots of the error $\|\betahat- \betastar\|_2$ after running
projected gradient descent on the nonconvex objective, with
sparsity $k \approx\sqrt{p}$. In all cases, we plotted the error
versus the rescaled sample size $\frac{n}{k \log p}$. As predicted
by Theorems \protect\ref{ThmStat} and \protect\ref{ThmOpt}, the
curves align for
different values of $p$ when plotted in this rescaled manner.
\textup{(a)}
Missing data case with i.i.d. covariates. \textup{(b)}~Vector autoregressive
data with additive noise. Each point represents an average over
100 trials.}
\label{SingleRegression}
\end{figure}

%
%co5 #&#
\begin{cors}
\label{CorGraphModel}
Suppose the columns of the matrix $\Theta$ are $k$-sparse, and suppose
the condition number $\kappa(\Theta)$ is nonzero and finite. Suppose
we have
%
%e3.12 #&#
\begin{equation}
\label{EqnDev2}
\bigl\|\widehat{\gamma}{}^{(j)} - \GamHat^{(j)}
\theta^j \bigr\|_\infty\leq \varphi(\Qprob, \sigeps) \sqrt{
\frac{\log\pdim}{\numobs}}\qquad \forall j,
\end{equation}
and suppose we have the following additional deviation condition on
$\SigHat$:
%
%e3.13 #&#
\begin{equation}
\label{ErrMaxBd} \|\SigHat- \Sigma\|_{\max} \le c\varphi(\Qprob,
\sigeps) \sqrt{\frac{\log p}{n}}.
\end{equation}
Finally, suppose the lower-RE condition holds uniformly over the
matrices~$\GamHat^{(j)}$\vadjust{\goodbreak} with the
scaling (\ref{EqnRSCConstants}). Then under the estimation procedure
of Algorithm~\ref{AlgMatrix}, there exists a universal constant $c_0$
such that
\[
\Vvert\Thetahat- \Theta\Vvert_{\mathrm{op}} \le \frac{c_0\kappa^2(\Sigma)}{\lammin(\Sigma)} \biggl(
\frac{\varphi
(\Qprob,
\sigeps)}{\lammin(\Sigma)} + \frac{\varphi(\Qprob,
\sigeps)}{\allow} \biggr) k\sqrt{\frac{\log p}{n}}.
\]
\end{cors}

\begin{rem*}
Note that Corollary~\ref{CorGraphModel} is again a deterministic
result, with parallel structure to Theorem~\ref{ThmStat}. Furthermore,
the deviation bounds (\ref{EqnDev2}) and (\ref{ErrMaxBd}) hold for all
scenarios considered in Section~\ref{SecConseq} above, using
Corollaries~\ref{CorAddNoise}--\ref{CorARmissing} for the first two inequalities, and a similar bounding
technique for $\|\SigHat- \Sigma\|_{\max}$; and the lower-RE
condition holds over all matrices $\GamHat^{(j)}$ by the same
technique used to establish the lower-RE condition for $\GamHat$. The
uniformity of the lower-RE bound over all sub-matrices holds because
\[
0 < \lambda_{\min}(\Sigma) \le\lambda_{\min}(
\Sigma_{-j,-j}) \le \lambda_{\max}(\Sigma_{-j,-j}) \le
\lambda_{\max}(\Sigma) < \infty.
\]
Hence, the error bound in Corollary~\ref{CorGraphModel} holds
with probability at least $1 - c_1 \exp(-c_2 \log p)$ when $n \succsim
\kdim\log\pdim$, for the appropriate values of $\varphi$ and
$\allow$.
\end{rem*}

%%%%%%%%%%%%%%%%%%%%%%%%%%%%%%%%%%%%%%%%%%%%%%%%%%%%%%%%%%%%%%%%%%%%%%%%

%s4 #&#
\section{Simulations}
\label{SecSims}

In this section, we report some additional simulation results to
confirm that the scalings predicted by our theory are sharp. In
Figure~\ref{FigRescaling} following Theorem~\ref{ThmStat}, we showed
that the error curves align when plotted against a suitably rescaled
sample size, in the case of additive noise perturbations. Panel (a)
of Figure~\ref{SingleRegression} shows these same types of rescaled
curves for the case of missing data, with sparsity $k \approx
\sqrt{\pdim}$, covariate matrix $\Sigma_x = I$, and missing
fraction\vadjust{\goodbreak}
$\rho= 0.2$, whereas panel (b) shows the rescaled plots for the
vector autoregressive case with additive noise perturbations, using
a driving matrix $A$ with $\Vvert A \Vvert_{\mathrm{op}} = 0.2$. Each
point corresponds
to an average over 100 trials. Once again, we see excellent agreement
with the scaling law provided by Theorem~\ref{ThmStat}.

\begin{figure}[b]
\begin{tabular}{@{}cc@{}}

\includegraphics{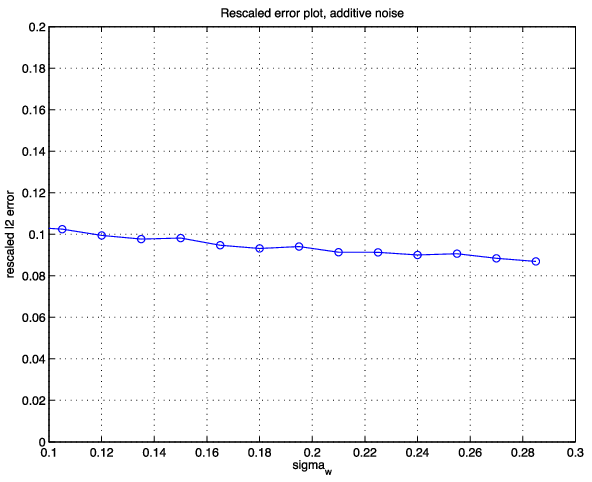}
 & \includegraphics{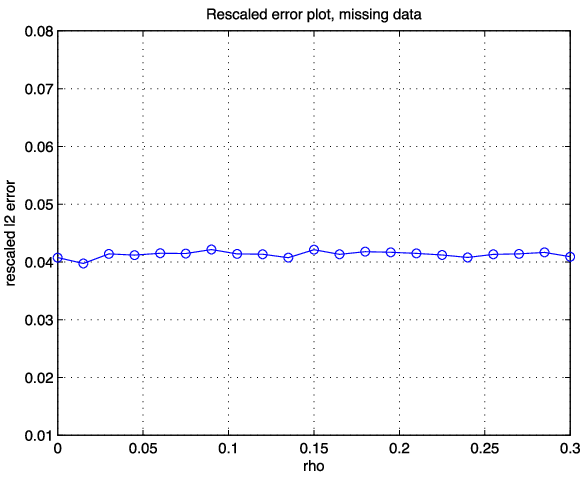}\\
(a) & (b)
\end{tabular}
\caption{\textup{(a)} Plot of the rescaled $\ell_2$-error
$\frac{\sqrt{1+\sigma_w^2}}{\sigma_w + 0.5} \|\betahat- \betastar
\|_2$
versus the additive noise standard deviation $\sigma_w$ for the
i.i.d. model with additive noise. \textup{(b)} Plot of the rescaled
$\ell_2$-error $\frac{\|\betahat- \betastar\|_2}{1 +
0.5(1-\rhoprob)}$ versus the missing fraction $\rhoprob$ for the
i.i.d. model with missing data. Both curves are roughly constant,
showing that our error bounds on $\|\betahat- \betastar\|_2$
exhibit the proper scaling. Each point represents an average over
200 trials.}
\label{ScalingPlot}
\end{figure}

We also ran simulations to verify the form of the function
$\varphi(\Qprob,\sigma_\varepsilon)$ appearing in
Corollaries~\ref{CorAddNoise} and~\ref{CorMAR}. In the additive noise
setting for i.i.d. data, we set $\Sigma_x = I$ and $\varepsilon$ equal
to i.i.d. Gaussian noise with $\sigma_\varepsilon= 0.5$. For a~fixed
value of the parameters $p = 256$ and $k \approx\log p$, we ran the
projected gradient descent algorithm for different values of $\sigma_w
\in(0.1,0.3)$, such that $\Sigma_w = \sigma_w^2 I$ and $n \approx60
(1+\sigma_w^2)^2 k \log p$, with $\|\betastar\|_2 = 1$. According to
the theory, $\frac{\varphi(\Qprob, \sigeps)}{\alpha}
\asymp
(\sigma_w + 0.5) \sqrt{1 + \sigma_w^2}$,
%1+\sigma_w^2 + \sigeps\sqrt{1 + \sigma_w^2}$,}
so that
\[
\bigl\|\betahat- \betastar\bigr\|_2 \precsim(\sigma_w + 0.5)
\sqrt{1 + \sigma_w^2} \sqrt{
\frac{k \log p}{(1+\sigma_w^2)^2 k \log p}} \asymp\frac{\sigma_w + 0.5}{\sqrt{1+\sigma_w^2}}.
\]
%
%%
%p}{(1+\sigma_w^2)^2 k \log p}} \\
%%
%& = 1 + \frac{0.5}{\sqrt{1+\sigma_w^2}}.
%%
%}
In order to verify this theoretical prediction, we
plotted $\sigma_w$ versus the rescaled error
$\frac{\sqrt{1+\sigma_w^2}}{\sigma_w + 0.5} \|\betahat-
\betastar\|_2$. %\pldel{\mbox{$\frac{\|\betahat- \betastar\|_2}{1
%+ \frac{0.5}{\sqrt{1 + \sigma_w^2}}}$.}}
As shown by\vspace*{2pt}
Figure~\ref{ScalingPlot}(a), the curve is roughly constant, as
predicted by the theory.\vspace*{1pt}

%f4 #&#

Similarly, in the missing data setting for i.i.d. data, we set
$\Sigma_x = I$ and~$\varepsilon$ equal to i.i.d. Gaussian noise with
$\sigma_\varepsilon= 0.5$. For a fixed value of the parameters $p = 128$
and $k \approx\log p$, we ran simulations for different values of the
missing data parameter $\rhoprob\in(0, 0.3)$, such that $n \approx
\frac{60}{(1-\rhoprob)^4} k \log p$. According to the theory,\vadjust{\goodbreak}
$\frac{\varphi(\Qprob, \sigma_\varepsilon)}{\alpha} \asymp
\frac{\sigeps}{1-\rho} + \frac{1}{(1-\rho)^2}$. Consequently, with
our specified scalings of $(\numobs, \pdim, \kdim)$, we should expect
a bound of the form
\[
\bigl\|\betahat- \betastar\bigr\|_2 \precsim\frac{\varphi(\Qprob,
\sigma_\varepsilon)}{\alpha} \sqrt{
\frac{k \log p}{n}} \asymp 1 + 0.5(1-\rho).\vspace*{-2pt}
\]
The plot\vspace*{1pt} of $\rhoprob$ versus the rescaled error $\frac{\|\betahat-
\betastar\|_2}{1 + 0.5(1-\rhoprob)}$ is shown in
Figure~\ref{ScalingPlot}(b). The curve is again roughly constant,
agreeing with theoretical results.

Finally, we studied the behavior of the inverse covariance matrix
estimation algorithm on three types of Gaussian graphical models:
\begin{longlist}[(a)]
\item[(a)] \textit{Chain-structured graphs.} In this case, all nodes of
the graph are arranged in a linear chain. Hence, each node (except
the two end nodes) has degree $k = 2$. The diagonal entries of
$\Theta$ are set equal to 1, and all entries corresponding to links
in the chain are set equal to $0.1$. Then $\Theta$ is rescaled so
$\Vvert\Theta\Vvert_{\mathrm{op}} = 1$.
\item[(b)] \textit{Star-structured graphs.} In this case, all nodes are
connected to a central node, which has degree $k \approx0.1p$. All
other nodes have degree 1. The diagonal entries of $\Theta$ are set
equal to 1, and all entries corresponding to edges in the graph are
set equal to $0.1$. Then $\Theta$ is rescaled so $\Vvert\Theta
\Vvert_{\mathrm{op}} =
1$.
\item[(c)] \textit{Erd\H{o}s--Renyi graphs.} This example comes from
Rothman et al.~\cite{Rot08}. For a sparsity parameter $k \approx
\log p$, we randomly generate the matrix $\Theta$ by first
generating the matrix $B$ such that the diagonal entries are 0, and
all other entries are independently equal to 0.5 with probability
$k/p$, and 0 otherwise. Then $\delta$ is chosen so that $\Theta= B
+ \delta I$ has condition number $p$. Finally, $\Theta$ is rescaled
so $\Vvert\Theta\Vvert_{\mathrm{op}} = 1$.
\end{longlist}
After generating the matrix $X$ of $n$ i.i.d. samples from the
appropriate graphical model, with covariance matrix $\Sigma_x =
\Theta^{-1}$, we generated the corrupted matrix $Z = X + W$ with
$\Sigma_w = (0.2)^2 I$ in the additive noise case, or the missing data
matrix $Z$ with $\rhoprob= 0.2$ in the missing data case.

%f5 #&#
\begin{figure}
\begin{tabular}{@{}cc@{}}

\includegraphics{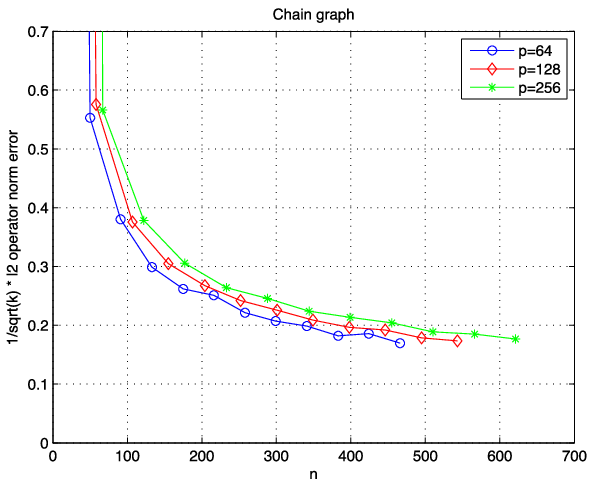}
 & \includegraphics{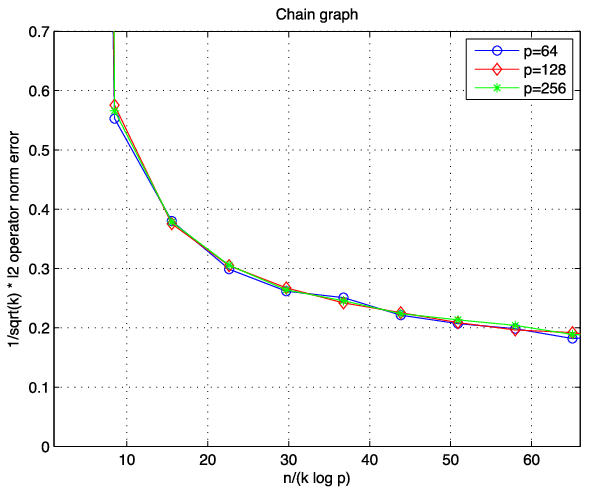}\\
(a) & (b)\\[6pt]

\includegraphics{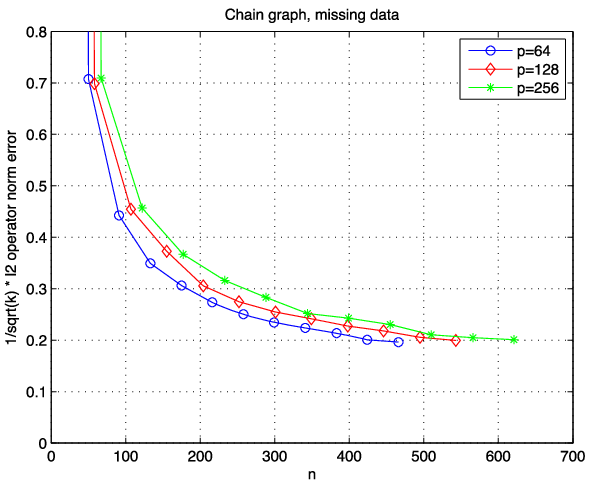}
 & \includegraphics{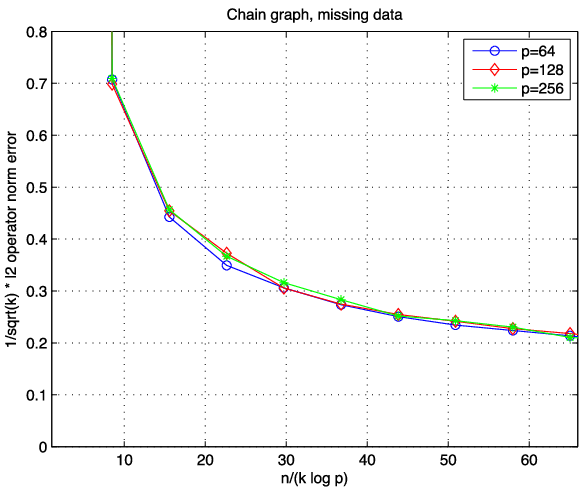}\\
(c) & (d)
\end{tabular}
\caption{\textup{(a)} Plots of the error $\Vvert\Thetahat- \Theta\Vvert_{\mathrm{op}}$ after
running projected gradient descent on the nonconvex objective for a
chain-structured Gaussian graphical model with additive noise. As
predicted by Theorems \protect\ref{ThmStat} and \protect\ref
{ThmOpt}, all curves
align when the error is rescaled by $\frac{1}{\sqrt{k}}$ and plotted
against the ratio $\frac{\numobs}{\kdim\log\pdim}$, as shown in
\textup{(b)}. Plots \textup{(c)} and \textup{(d)} show the results of simulations on missing
data sets. Each point represents the average over 50 trials.}
\label{ChainPlots}
\end{figure}

Panels (a) and (c) in Figure~\ref{ChainPlots} show the rescaled
$\ell_2$-error $\frac{1}{\sqrt{k}}\Vvert\Thetahat- \Theta\Vvert_{\mathrm{op}}$
plotted against the sample size $n$ for a chain-structured graph. In
panels (b) and (d), we have $\ell_2$-error plotted against the
rescaled sample size, $n/(k \log p)$. Once again, we see good
agreement with the theoretical predictions. We have obtained
qualitatively similar results for the star and Erd\H{o}s--Renyi
graphs.\vspace*{-2pt}

%%%%%%%%%%%%%%%%%%%%%%%%%%%%%%%%%%%%%%%%%%%%%%%%%%%%%%%%%%%%%%%%%%%%%%%%%%%%%%

%s5 #&#
\section{Proofs}
\label{SecProofs}

In this section, we prove our two main theorems. For the more
technical proofs of the corollaries, see the supplementary Appendix~\cite{LohWai11}.\vspace*{-2pt}

%s5.1 #&#
\subsection{\texorpdfstring{Proof of Theorem \protect\ref{ThmStat}}{Proof of Theorem 1}}

Let $\Loss(\beta) = \frac{1}{2} \beta^T \Gammahat\beta-
\langle\widehat{\gamma}, \beta
\rangle+ \regpar\|\beta\|_1$ denote the loss
function to be minimized. This definition captures both the
estimator~(\ref{EqnGeneral}) with $\regpar= 0$ and the
estimator~(\ref{EqnGeneralSlight}) with the choice of $\regpar$ given
in the theorem statement. For either estimator, we are guaranteed
that $\betastar$ is feasible and $\betahat$ is optimal for the
program, so $\Loss(\betahat) \le\Loss(\betastar)$. Indeed,\vspace*{1pt} in the
regularized case, the $\kdim$-sparsity of $\betastar$ implies that
$\|\betastar\|_1 \le\sqrt{\kdim} \|\betastar\|_2 \le\ANNOYCON
\sqrt{\kdim}$. Defining the error vector\vadjust{\goodbreak} $\nuhat\defn\betahat-
\betastar$ and performing some algebra leads to the equivalent
inequality
%
%e5.1 #&#
\begin{equation}
\label{EqnBasic} \tfrac{1}{2} \nuhat^T \Gammahat\nuhat\le
\bigl\langle \nuhat, \widehat{\gamma} - \Gammahat\betastar \bigr\rangle+ \regpar \bigl\{ \bigl\|
\betastar\bigr\|_1 - \bigl\|\betastar+ \nuhat\bigr\|_1 \bigr\}.
\end{equation}
In the remainder of the proof, we first derive an upper bound for
the right-hand side of this inequality. We then use this upper bound
and the
lower-RE condition to show that the error vector $\nuhat$ must satisfy
the inequality
%
%e5.2 #&#
\begin{equation}
\label{EqnConeImply} \|\nuhat\|_1 \leq8 \sqrt{\kdim} \|\nuhat
\|_2.
\end{equation}
Finally, we combine inequality (\ref{EqnConeImply}) with the
lower-RE condition to derive a lower bound on the left-hand side of the basic
inequality (\ref{EqnBasic}). Combined with our earlier upper bound on
the right-hand side, some algebra yields the claim.\vadjust{\goodbreak}

\subsubsection*{Upper bound on right-hand side} We first upper-bound the
right-hand side of
inequality (\ref{EqnBasic}). H\"{o}lder's inequality gives
$\langle\nuhat, \widehat{\gamma}-
\Gammahat\betastar\rangle\le\|\nuhat\|_1
\|\widehat{\gamma}- \Gammahat\betastar\|_\infty$. By
the triangle
inequality, we have
\[
\bigl\|\widehat{\gamma}- \Gammahat\betastar\bigr\|_\infty\leq
\bigl\|\widehat{\gamma}-
\Sigma_x \betastar\bigr\|_\infty+ \bigl\|(\Sigma_x -
\Gammahat) \betastar\bigr\|_\infty\stackrel{(\mathrm{i})} {\leq} 2 \varphi(\Qprob,
\sigeps) \sqrt{\frac{\log p}{n}},
\]
where inequality (i) follows from the deviation
conditions (\ref{EqnDev}). Combining the pieces, we conclude that
%
%e5.3 #&#
\begin{eqnarray}
\label{EqnBasicUpper}
\bigl\langle\nuhat, \widehat{\gamma}- \Gammahat\betastar\bigr\rangle
&\le&2 \|\nuhat\|_1 \varphi(\Qprob, \sigeps)
\sqrt{\frac{\log\pdim}{\numobs}}\nonumber\\[-8pt]\\[-8pt]
&=& \bigl( \|\nuhat_S\|_1 + \|\nuhat_{\Sbar}
\|_1 \bigr) 2 \varphi (\Qprob, \sigeps)
\sqrt{\frac{\log\pdim}{\numobs}}.\nonumber
\end{eqnarray}
On the other hand, we have
%
%e5.4 #&#
\begin{eqnarray}
\label{EqnLowerReg}
\bigl\|\betastar+ \nuhat\bigr\|_1 - \bigl\|\betastar
\bigr\|_1 &\geq& \bigl\{ \bigl\|\betastar_S\bigr\|_1 - \|
\nuhat_S\|_1 \bigr\} + \|\nuhat_{\Sbar}
\|_1 - \bigl\|\betastar\bigr\|_1 \nonumber\\[-8pt]\\[-8pt]
&=& \|\nuhat_{\Sbar}
\|_1 - \|\nuhat_S\|_1,\nonumber
\end{eqnarray}
where we have exploited the sparsity of $\betastar$ and applied the
triangle inequality. Combining the pieces, we conclude that the
right-hand side
of inequality~(\ref{EqnBasic}) is upper-bounded by
%
%e5.5 #&#
\begin{equation}
\label{Eqnright-handsideUpper} 2 \varphi(\Qprob, \sigeps) \sqrt{
\frac{\log\pdim}{\numobs}} \bigl(\|\nuhat_S\|_1 + \|
\nuhat_{S^c}\|_1\bigr) + \regpar\bigl\{\|\nuhat_S
\|_1 - \|\nuhat_{S^c}\|_1\bigr\},
\end{equation}
a bound that holds for any nonnegative choice of $\regpar$.

\subsubsection*{\texorpdfstring{Proof of inequality (\protect\ref{EqnConeImply})}
{Proof of inequality (5.2)}}

We first consider the constrained program~(\ref{EqnGeneral}), with $R
= \|\betastar\|_1$, so $\|\betahat\|_1 = \|\betastar+ \nuhat\|_1
\leq\|\betastar\|_1$. Combined with inequality~(\ref{EqnLowerReg}),
we conclude that $\|\nuhat_{\Sbar}\|_1 \leq\|\nuhat_{S}\|_1$.
Consequently, we have the inequality $\|\nuhat\|_1 \leq2
\|\nuhat_S\|_1 \leq2 \sqrt{\kdim} \|\nuhat\|_2$, which is a
slightly stronger form of the bound (\ref{EqnConeImply}).

For the regularized estimator (\ref{EqnGeneralSlight}), we first note
that our choice of $\regpar$ guarantees that the
term (\ref{Eqnright-handsideUpper}) is at most $\frac{3 \regpar}{2}
\|\nuhat_S\|_1 - \frac{\regpar}{2} \|\nuhat_{S^c}\|_1$. Returning to
the basic inequality, we apply the lower-RE condition to
lower-bound the left-hand side, thereby obtaining the inequality
\[
- \frac{\taulow}{2} \|\nuhat\|_1^2 \leq
\frac{1}{2} \bigl(\allow \|\nuhat\|_2^2 - \taulow \|
\nuhat\|_1^2 \bigr) \leq\frac{3 \regpar}{2} \|
\nuhat_S\|_1 - \frac{\regpar}{2} \|
\nuhat_{S^c}\|_1.
\]
By the triangle inequality, we have $\|\nuhat\|_1 \leq\|\betahat\|_1
+ \|\betastar\|_1 \leq2 \ANNOYCON\sqrt{\kdim}$. Since we have
assumed $\sqrt{\kdim} \taulow(\numobs, \pdim) \leq
\frac{\varphi(\Qprob, \sigeps)}{\ANNOYCON} \sqrt{\frac{\log
\pdim}{\numobs}} $, we are guaranteed that
\[
\frac{\taulow(\numobs, \pdim)}{2} \|\nuhat\|_1^2 \leq\varphi (\Qprob,
\sigeps) \sqrt{\frac{\log\pdim}{\numobs}} \|\nuhat\|_1 \leq
\frac{\regpar}{4} \|\nuhat\|_1
\]
by our choice of $\regpar$. Combining the pieces, we conclude
that
\[
0 \leq\frac{3 \regpar}{2} \|\nuhat_S\|_1 -
\frac{\regpar}{2} \|\nuhat_{\Sbar}\|_1 + \frac{\regpar}{4} \bigl(
\|\nuhat_S\|_1 + \|\nuhat_{\Sbar}
\|_1 \bigr) = \frac{7 \regpar}{4} \|\nuhat_S\|_1
- \frac{\regpar}{4} \|\nuhat_{\Sbar}\|_1
\]
and rearranging implies $\|\nuhat_{\Sbar}\|_1 \leq7 \|\nuhat_S\|_1$,
from which we conclude that $\|\nuhat\|_1 \leq8 \sqrt{\kdim}
\|\nuhat\|_2$, as claimed.

\subsubsection*{Lower bound on left-hand side}
We now derive a lower bound on the left-hand side of
inequality (\ref{EqnBasic}). Combining inequality (\ref{EqnConeImply})
with the RE condition (\ref{EqnRSC}) gives
%
%e5.6 #&#
\begin{equation}
\label{EqnBasicLower}\qquad
\nuhat^T \Gammahat\nuhat\ge\allow\|\nuhat
\|_2^2 - \taulow (\numobs, \pdim) \|\nuhat
\|_1^2 \geq \bigl\{ \allow- 64 \kdim \taulow(\numobs,
\pdim) \bigr\} \|\nuhat\|_2^2 \geq \frac{\allow}{2} \|
\nuhat\|_2^2,
\end{equation}
where the final step uses our assumption that $\kdim\taulow(\numobs,
\pdim) \leq\frac{\allow}{128}$.

Finally, combining bounds (\ref{Eqnright-handsideUpper}), (\ref
{EqnConeImply})
and (\ref{EqnBasicLower}) yields
\begin{eqnarray*}
\frac{\allow}{4} \|\nuhat\|_2^2 & \le & 2 \max \Biggl\{2
\varphi (\Qprob, \sigeps) \sqrt{\frac{\log p}{n}}, \regpar \Biggr\} \|\nuhat
\|_1
\\
& \le & 32 \sqrt{k} \max \Biggl\{\varphi(\Qprob, \sigeps) \sqrt{
\frac{\log p}{n}}, \regpar \Biggr\} \|\nuhat\|_2,
\end{eqnarray*}
giving inequality (\ref{EqnEllTwoBound}). Using
inequality (\ref{EqnConeImply}) again gives
inequality (\ref{EqnEllOneBound}).

%%%%%%%%%%%%%%%%%%%%%%%%%%%%%%%%%%%%%%%%%%%%%%%%%%%%%%%%%%%%%%%%%%%%%%%%%%%%%%

%s5.2 #&#
\subsection{\texorpdfstring{Proof of Theorem \protect\ref{ThmOpt}}{Proof of Theorem 2}}

We begin by proving the claims for the constrained problem, and
projected gradient descent. For the $\ell_2$-error bound, we make use
of Theorem 1 in the pre-print of Agarwal et al.~\cite{AgaEtal11}.
Their theory, as originally stated, requires that the loss function be
convex, but a careful examination of their proof shows that
their arguments hinge on restricted strong convexity and
smoothness assumptions, corresponding to a more general version of
the lower- and upper-RE conditions given here. Apart from these
conditions, the proof exploits the fact that the sub-problems
defining the gradient updates (\ref{EqnProjGrad})
and (\ref{EqnMirror}) are convex. Since the loss function itself
appears only in a linear term, their theory still applies.
%is convex.}

In order to apply Theorem 1 in their paper, we first need to compute
the tolerance parameter $\varepsilon^2$ defined there; since $\betastar$
is supported on the set $S$ with $|S| = \kdim$ and the RE conditions
hold with $\tau\asymp\frac{\log\pdim}{\numobs}$, we find that
\begin{eqnarray*}
\varepsilon^2 & \leq & \plaincon\frac{\log
\pdim}{\alup\numobs} \bigl( \sqrt{\kdim} \bigl\|
\betahat- \betastar\bigr\|_2 + 2 \bigl\|\betahat- \betastar\bigr\|_1
\bigr)^2
\\
& \leq & \plaincon_2' \frac{\kdim\log\pdim}{\alup
\numobs} \bigl\|
\betahat- \betastar\bigr\|_2^2 + \plaincon_1
\frac{\log
\pdim}{\alup\numobs} \bigl\|\betahat- \betastar\bigr\|_1^2
\\
& \leq & \plaincon_2 \bigl\|\betahat- \betastar
\bigr\|_2^2 + \plaincon_1 \frac{\log\pdim}{\alup\numobs} \bigl\|
\betahat- \betastar\bigr\|_1^2,
\end{eqnarray*}
where the final inequality makes use of the assumption that $\numobs
\succsim\kdim\log\pdim$. Similarly, we may compute the
contraction coefficient to be
%
%e5.7 #&#
\begin{equation}
\label{EqnMyContract} \gamma= \biggl(1 - \frac{\allow}{\alup} + \frac{c_1 k\log p}{\alup
n}
\biggr) \biggl(1 - \frac{c_2 k\log p}{\alup n} \biggr)^{-1},
\end{equation}
so $\gamma\in(0,1)$ for $n \succsim\kdim\log\pdim$.

We now establish the $\ell_1$-error bound. First, let $\Delta^t \defn
\beta^t - \betastar$. Since $\beta^t$ is feasible and $\betahat$ is
optimal with an active constraint, we have $\|\beta^t\|_1 \leq
\|\betahat\|_1$. Applying the triangle inequality gives
\begin{eqnarray*}
\|\betahat\|_1 & \le & \bigl\|\betastar\bigr\|_1 + \bigl\|\betahat-
\betastar\bigr\|_1 = \bigl\|\betastar_S\bigr\|_1 + \bigl\|
\betahat- \betastar\bigr\|_1,
\\
\bigl\|\beta^t\bigr\|_1 &=& \bigl\|\betastar+
\Delta^t\bigr\|_1 \ge \bigl\|\betastar_S +
\Delta^t_{S^c}\bigr\|_1 - \bigl\|\Delta^t_S
\bigr\|_1 = \bigl\|\betastar_S\bigr\|_1 + \bigl\|
\Delta^t_{S^c}\bigr\|_1 - \bigl\|\Delta^t_S
\bigr\|_1;
\end{eqnarray*}
combining the bounds yields $\|\Delta^t_{S^c}\|_1 \le\|\Delta^t_S\|_1
+ \|\betahat- \betastar\|_1$. Then
\[
\bigl\|\Delta^t\bigr\|_1 \le2\bigl\|\Delta^t_S
\bigr\|_1 + \bigl\|\betahat- \betastar\bigr\|_1 \le 2\sqrt{\kdim} \bigl\|
\Delta^t\bigr\|_2 + \bigl\|\betahat- \betastar\bigr\|_1,
\]
so
\[
\bigl\|\beta^t - \betahat\bigr\|_1 \le\bigl\|\betahat- \betastar
\bigr\|_1 + \bigl\|\Delta^t\bigr\|_1 \le2\sqrt{\kdim}
\bigl(\bigl\|\beta^t - \betahat\bigr\|_2 + \bigl\|\betahat- \betastar
\bigr\|_2\bigr) + 2\bigl\|\betahat- \betastar\bigr\|_1.
\]

Turning to the Lagrangian version, we exploit Theorem 2 in Agarwal et
al.~\cite{AgaEtal11}, with $\mathcal{M}$ corresponding to the subspace of all vectors
with support contained within the support set of $\betastar$. With
this choice, we have $\psi(\mathcal{M}) = \sqrt{k}$, and the
contraction coefficient $\CONTRAC$ takes the previous
form (\ref{EqnMyContract}), so that the assumption $n \succsim\kdim
\log\pdim$ guarantees that $\CONTRAC\in(0,1)$. It remains to
verify that the requirements are satisfied. From the conditions in
our Theorem~\ref{ThmOpt} and using the notation of Agarwal et al.
\cite{AgaEtal11}, we have
$\beta(\mathcal{M}) = \order(\frac{\log\pdim}{\numobs})$ and
$\overline{\rho} = \sqrt{k}$, and the condition $\numobs\succsim
\kdim
\log\pdim$ implies that $\xi(\mathcal{M}) = \order(1)$.
Putting together the pieces, we find that the compound tolerance
parameter $\varepsilon^2$ satisfies the bound $\varepsilon^2 = \order(
\frac{\kdim\log\pdim}{\numobs} \|\betahat- \betastar\|_2^2) =
\order(\|\betahat- \betastar\|_2^2)$, so the claim follows.

%s6 #&#
\section{Discussion}

In this paper, we formulated an $\ell_1$-constrained minimization
problem for sparse linear regression on corrupted data. The source of
corruption may be additive noise or missing data, and although the
resulting objective is not generally convex, we showed that projected
gradient descent is guaranteed to converge to a point within
statistical precision of the optimum. In addition, we established
$\ell_1$- and $\ell_2$-error bounds that hold with high probability
when the data are drawn i.i.d. from a sub-Gaussian distribution, or
drawn from a Gaussian vector autoregressive process. Finally, we
applied our methods to sparse inverse covariance estimation for a
Gaussian graphical model with corruptions, and obtained spectral norm
rates of the same order as existing rates for uncorrupted, i.i.d.
data.

Future directions of research include studying more general types of
dependencies or corruption in the covariates of regression, such as
more general types of multiplicative noise, and performing sparse
linear regression for corrupted data with additive noise when the
noise covariance is unknown and replicates of the data may be
unavailable. As pointed out by a reviewer, it would also be
interesting to study the performance of our algorithms on data that
are not sub-Gaussian, or even under model mismatch. In addition,
one might consider other loss functions, where it is more difficult to
correct the objective for corrupted covariates. Finally, it remains
to be seen whether or not our techniques---used to show that certain
nonconvex problems can solved to statistical precision---can be
applied more broadly.

\section*{Acknowledgments}

The authors thank Alekh Agarwal, Sahand Negahban, John Duchi and Alexandre Tsybakov for
useful discussions and guidance. They are also grateful to the
Associate Editor and anonymous referees for improvements on the paper.

\begin{supplement}%[id=suppA]
\stitle{Supplementary material for:
High-dimensional regression with noisy and
missing data: Provable guarantees with nonconvexity\\}
\slink[doi]{10.1214/12-AOS1018SUPP} %[doi,text={...}] - jei reikia
%suskaldyti doi
\sdatatype{.pdf}
\sfilename{aos1018\_supp.pdf}
\sdescription{Due to space constraints, we have relegated technical details of the
remaining proofs to the supplement~\cite{LohWai11}.}
\end{supplement}

% imsref loaded by lrinkeviciute, 2012-07-04 09:59:18
% imsref loaded by lrinkeviciute, 2012-07-04 10:25:27

\printaddresses

\end{document}